\pdfoutput=1
\RequirePackage{ifpdf}
\ifpdf 
\documentclass[pdftex]{sigma}
\else
\documentclass{sigma}
\fi

\numberwithin{equation}{section}

\newtheorem{Theorem}{Theorem}[section]
\newtheorem{Corollary}[Theorem]{Corollary}
\newtheorem{Lemma}[Theorem]{Lemma}
\newtheorem{Proposition}[Theorem]{Proposition}
{ \theoremstyle{definition}

\newtheorem{Note}[Theorem]{Note}

 }

\usepackage{enumerate}

\begin{document}
\allowdisplaybreaks

\newcommand{\arXivNumber}{1803.11230}

\renewcommand{\thefootnote}{}

\renewcommand{\PaperNumber}{095}

\FirstPageHeading

\ShortArticleName{Tronqu{\'e}e Solutions of the Third and Fourth Painlev{\'e} Equations}

\ArticleName{Tronqu{\'e}e Solutions\\ of the Third and Fourth Painlev{\'e} Equations\footnote{This paper is a~contribution to the Special Issue on Painlev\'e Equations and Applications in Memory of Andrei Kapaev. The full collection is available at \href{https://www.emis.de/journals/SIGMA/Kapaev.html}{https://www.emis.de/journals/SIGMA/Kapaev.html}}}

\Author{Xiaoyue XIA}

\AuthorNameForHeading{X.~Xia}

\Address{Department of Mathematics, The Ohio State University,\\ 100 Math Tower, 231 West 18th Avenue, Columbus OH, 43210-1174, USA}
\Email{\href{mailto:xiaxiaoyue9@gmail.com}{xiaxiaoyue9@gmail.com}}

\ArticleDates{Received April 04, 2018, in final form August 30, 2018; Published online September 08, 2018}

\Abstract{Recently in a paper by Lin, Dai and Tibboel, it was shown that the third and fourth Painlev{\'e} equations have tronqu{\'e}e and tritronqu{\'e}e solutions. We obtain global information about these tronqu{\'e}e and tritronqu{\'e}e solutions. We find their sectors of analyticity, their Borel summed representations in these sectors as well as the asymptotic position of the singularities near the boundaries of the analyticity sectors. We also correct slight errors in the paper mentioned.}

\Keywords{the third and fourth Painlev{\'e} equations; asymptotic position of singularities; tronqu{\'e}e solutions; tritronqu{\'e}e solutions; Borel summed representation}

\Classification{34M25; 34M40; 34M55}

\renewcommand{\thefootnote}{\arabic{footnote}}
\setcounter{footnote}{0}

\section{Introduction}\label{intro}
The well-known Painlev{\'e} equations were first introduced by Painlev{\'e} more than a century ago and have been investigated by many researchers. The Painlev{\'e} equations define new functions called Painlev{\'e} transcendents, which are considered as special nonlinear functions and their asymptotic behavior is of particular importance. For an overview of Painlev{\'e} equations and the asymptotic behavior of Painlev{\'e} transcendents please see, e.g.,~\cite{Clarkson} and~\cite{FIKN}. In recent decades there has been revived interest in Painlev{\'e} equations as they play important roles in various mathematical and physical applications (see, e.g., \cite{Clarkson,CC2,FFW,FIKN,K2,K3} for references to applications).

Boutroux first studied a family of particular solutions of the first Painlev{\'e} equations ${\rm P_I}$, which he named ``tronqu{\'e}e'' and ``tritronqu{\'e}e'' solutions in~\cite{Boutroux}. These special solutions of Painlev{\'e} equations have pole-free sectors while generic solutions have poles accumulating at $\infty$ in all sectors. Tronqu{\'e}e and tritronqu{\'e}e solutions receive attention not only for their interesting analytic property but also because they appear in a number of problems such as the Ising model~\cite{MTW}, the critical behavior in the NLS/Toda lattices~\cite{Dubrovin, DGK} and the analysis of the cubic oscillator~\cite{Masoero}.

For the first Painlev{\'e} equation some pioneering works based on the powerful techniques of isomonodromic deformation and reduction to Riemann--Hilbert problem were done in the study of tronqu{\'e}e solutions by Kapaev and coauthors. In \cite{Kapaev} and \cite{KK} the Stokes constant for the tritronqu{\'e}e solution of ${\rm P_I}$ was calculated for the first time. In \cite{K2} the global asymptotic behavior of the tronqu{\'e}e solutions of ${\rm P_I}$ was described with connection formulae presented. In~\cite{K3} and~\cite{IK} the global asymptotic behavior of the tronqu{\'e}e solutions of $\rm P_{II}$ was described with connection formulae presented. In \cite{K4} the global asymptotics of the solutions of the fourth Painlev{\'e} equation $\rm P_{IV}$ including its tronqu{\'e}e solutions was analyzed in detail. In \cite{DK} an fourth-order nonlinear ODE which controls the pole dynamics in the general solution of equation $\mathrm{P}_{\mathrm{I}}^{2}$ was studied. See also the monograph~\cite{FIKN} for a summary of recent developments in the theory of Painlev{\'e} equations based on this Riemann--Hilbert-isomonodromy method.

There is an impressive body of work on tronqu{\'e}e solutions and we only mention a few contributions here. Using approaches different from the Riemann--Hilbert-isomonodromy method Costin and coauthors analyzed tronqu{\'e}e solutions of ${\rm P_I}$ in \cite{CHT} and \cite{CCH} and obtained similar results to those in \cite{Kapaev} and \cite{KK}. In~\cite{GKK} the existence of the tritronqu{\'e}e solutions of~$\mathrm{P}_{\mathrm{I}}^{2}$, the second member in the $\rm P_I$ hierarchy was proved. In~\cite{Mazzocco} the existence of tronqu{\'e}e solutions of the second Painlev{\'e} hierarchy was proved. For the location of poles for the Hasting--McLeod solution to the second Painlev{\'e} equation please see~\cite{HXZ}, in which a special case of Novokshenov conjecture~\cite{Novo} was also proved. For the tronqu{\'e}e solutions to the third Painlev{\'e} equation please see~\cite{Lin}, which followed the idea in~\cite{Joshi}.

In this paper the tronqu{\'e}e and tritronqu{\'e}e solutions of the third and fourth Painlev{\'e} equation are studied:
\begin{gather} 
{\rm P_{III}}\colon \ {\frac {{\rm d}^{2}y}{{\rm d}{x}^{2}}} = \frac{1}{y}\left(\frac{{\rm d}y}{{\rm d}x}\right)^2-\frac{1}{x}\frac{{\rm d}y}{{\rm d}x}+\frac{1}{x}\left(\alpha y^2+\beta \right)+\gamma y^3+\frac{\delta}{y},
\nonumber\\
\label{P4}
{\rm P_{IV}}\colon \ {\frac {{\rm d}^{2}y}{{\rm d}{x}^{2}}} = \frac{1}{2y}\left(\frac{{\rm d}y}{{\rm d}x}\right)^2+\frac{3}{2}y^3+4xy^2+2\left(x^2-\alpha \right)y+\frac{\beta}{y},
\end{gather}
where $\alpha$, $\beta$, $\gamma$ and $\delta$ are arbitrary complex numbers. By B{\"a}cklund transformations (see \cite{Milne}) ${\rm P_{III}}$ can be reduced to
\begin{gather}\label{P31}
{\rm P}^{(i)}_{\rm III}\colon \ {\frac {{\rm d}^{2}y}{{\rm d}{x}^{2}}} = \frac{1}{y}\left(\frac{{\rm d}y}{{\rm d}x}\right)^2-\frac{1}{x}\frac{{\rm d}y}{{\rm d}x}+\frac{1}{x}\big(\alpha y^2+\beta \big)+ y^3-\frac{1}{y},\\
\label{P32}
{\rm P}^{(ii)}_{\rm III}\colon \ {\frac {{\rm d}^{2}y}{{\rm d}{x}^{2}}} = \frac{1}{y}\left(\frac{{\rm d}y}{{\rm d}x}\right)^2-\frac{1}{x}\frac{{\rm d}y}{{\rm d}x}+\frac{1}{x}\big( y^2+\beta \big)-\frac{1}{y}.
\end{gather}

In a famous paper \cite{MTW} by McCoy, Tracy and Wu, a one-parameter family of tronqu{\'e}e solutions of a special case of~\eqref{P31} where $\alpha = 2\nu$, $\beta = -2\nu$ was constructed, whose asymptotics at $\infty$ was congruent to ours~(\eqref{transseriesh} and~\eqref{P31cov}) and asymptotic expansion for small~$x$ was obtained. Furthermore, in a recent paper~\cite{FFW} by Fasondini et al.\ a comprehensive computer simulation of the McCoy--Tracy--Wu solution was given. The computer pictures of the pole distributions in~\cite{FFW} provide a good illustration of our description of the asymptotic position of poles in, e.g.,~\eqref{poledistrn}.

\eqref{P32} was studied as the degenerate $\rm P_{III}$ in \cite{KV1} and \cite{KV2}, and the position of the first array of poles was found in~\cite{KV2} via isomonodromy methods.

We base our methods on the results in \cite{OCIMRN} and \cite{OCInv}, which used the technique of Borel summation to describe the Stokes phenomenon. We obtain representations of tronqu{\'e}e solutions as Borel summed transseries (see also~\cite{OCCONM}), as well as the position of the first array of poles, bordering the sector of analyticity. We will first use a simple example to briefly illustrate some concepts in the Borel summation method. Please see also \cite[Section~5]{CCH} for an introduction.

In the following, we denote by $\mathcal{L}_{\phi}$ the Laplace transform
\begin{gather*}
f\longmapsto \int_{0}^{\infty {\rm e}^{{\rm i}\phi}} f(p){\rm e}^{-xp}{\rm d}p,
\end{gather*}
where $\phi \in \mathbb{R}$. See also \cite[p.~8]{OCIMRN} for the notation.

Assume that we have a formal series
\begin{gather*}
\tilde{f}(w) = \sum_{n=0}^{\infty}a_n w^{-r-n},\qquad \operatorname{Re}(r)>0,
\end{gather*}
where the series $\sum\limits_{n=0}^{\infty}a_n x^{n}$ has a positive radius of convergence. The Borel transform of $\tilde{f}$ is defined to be the formal power series
\begin{gather*}
\big(\mathcal{B} \tilde{f} \big)(p) := \sum_{n=0}^{\infty}\frac{a_n p^{n+r-1}}{\Gamma(n+r)}.
\end{gather*}

In most cases the explicit solution of a differential equation is not known. We may obtain classical asymptotic series as formal power series solutions, but these formal solutions do not contain parameters that help us distinguish between actual solutions. This is illustrated in the following simple ordinary differential equation at the irregular singularity at $x=\infty$
\begin{gather}\label{linearexample}
y'+y=\frac{1}{x^2}.
\end{gather}
The unique formal power series solution for $x\to \infty$ is
\begin{gather*}
\tilde{y}_0(x) = \sum_{n=1}^{\infty}\frac{n!}{x^{n+1}}
\end{gather*}
and the general solution to \eqref{linearexample} is
\begin{gather*}
y(x;C) = y_0(x)+C{\rm e}^{-x}, \qquad \textrm{where} \quad y_0(x) = {\rm e}^{-x} \int_{x_1}^{x}\frac{{\rm e}^s}{s^2} {\rm d}s \sim \tilde{y}_0(x) \quad \textrm{as} \quad x\to +\infty.
\end{gather*}

The idea of transseries solution is a completion of classical formal power series solution in the sense that the transseries solution representation includes the free parameters which appear in the actual solutions. In the example above, if we let
\begin{gather}\label{transseries0}
\tilde{y}(x) = \tilde{y}_0(x) + C {\rm e}^{-x} \qquad \textrm{for} \quad x\to +\infty,
\end{gather}
then \eqref{transseries0} is a formal solution to \eqref{linearexample} and the simplest example of a transseries.

Under appropriate conditions (see \cite{OCIMRN,OCDuke,OCInv}), given $\phi$, the operator $\mathcal{L}_{\phi}\mathcal{B}$ is a one-to-one map between the transseries solutions and actual solutions. In the example~\eqref{linearexample}, the actual solutions have representation:
\begin{gather*}
y(x) =
\begin{cases}
\mathcal{L}_{\phi} \mathcal{B}\tilde{y}_0(x) + C_{+} {\rm e}^{-x}, & -\phi = \arg(x) \in \big(0,\tfrac{\pi}{2}\big),\\
\mathcal{L}_{\phi} \mathcal{B}\tilde{y}_0(x) + C_{-} {\rm e}^{-x}, & -\phi = \arg(x) \in \big({-}\frac{\pi}{2},0\big).
\end{cases}
\end{gather*}
The value $C_+-C_-$ is called the Stokes constant. This representation is a trivial example of Borel summed representation of solutions. In the case of nonlinear systems such as~\eqref{NF}, the transseries solution is of the form~\eqref{NFsys} and the Borel summed representation of actual solutions is of the form~\eqref{solnsys}.

In this paper we study tronqu{\'e}e solutions of~\eqref{P31}, \eqref{P32} and~\eqref{P4} by first transforming each of them into a~second-order differential equation of the following form
\begin{gather}\label{eqh}
h''(w) - h(w)+\frac{1}{w}\big[\left(\beta_2-\beta_1\right) h(w)+\left(\beta_2+\beta_1\right) h'(w)\big] = g(w,h,h'),
\end{gather}
where $\beta_1$ and $\beta_2$ are constants and $g(w,h,h')$ is analytic at $(\infty,0,0)$. See \eqref{P31cov} for the change of variable for \eqref{P31}, see \eqref{P32cov} for the change of variable for \eqref{P32} and see \eqref{P41cov}, \eqref{P42cov} and \eqref{P43cov} for the change of variable for~\eqref{P4}. Next we make the substitution
\begin{gather}\label{htou}
\begin{bmatrix}
h(w)\\
h'(w)
\end{bmatrix}=
\begin{bmatrix}
1-\dfrac{\beta_1}{2w}&1+\dfrac{\beta_2}{2w}\vspace{1mm}\\
-1-\dfrac{\beta_1}{2w}&1-\dfrac{\beta_2}{2w}
\end{bmatrix}
\mathbf{u}(w).
\end{gather}
Then $\mathbf{u}$ is a solution to the following normalized (see \cite{OCInv}) 2-dimensional differential system:
\begin{gather}\label{NF}
\mathbf{u}'+\left(\hat{\Lambda} +\frac{\hat{B}}{w}\right) \mathbf{u}=\mathbf{g} (w,\mathbf{u} ),
\end{gather}
where
\begin{gather*}
\hat{\Lambda} =
\begin{bmatrix}
1&0\\
0&-1
\end{bmatrix},\qquad
\hat{B} =
\begin{bmatrix}
\beta_1&0\\
0&\beta_2
\end{bmatrix},
\end{gather*}
and $\mathbf{g} (w,\mathbf{u} )$ is analytic at $(\infty,\mathbf{0})$ with $\mathbf{g}(w,\mathbf{u}) = O\big(w^{-2}\big)+O\big(|\mathbf{u}|^2\big)$ as $w\to \infty$ and $\mathbf{u} \to \mathbf{0}$.

We obtain information about the tronqu{\'e}e and tritronqu{\'e}e solutions of the normalized sys\-tem~\eqref{NF} such as their existence, regions of analyticity and asymptotic position of poles through which we obtain corresponding results regarding tronqu{\'e}e and tritronqu{\'e}e solutions of~$\rm P_{III}$ and~$\rm P_{IV}$. See also~\cite{CCH} in which a similar approach was used to study the tronqu{\'e}e solutions of the first Painlev{\'e} equation.

\section{Tronqu{\'e}e solutions of (\ref{eqh})}\label{NFinfo}

\subsection{Formal solutions and tronqu{\'e}e solutions of (\ref{eqh})}

In Proposition~\ref{NFformal} and Theorem~\ref{NFactual} we present formal and actual solutions on the right half $w$-plane $S_1 := \big\{w\colon \arg(w) \in \big({-}\frac{\pi}{2},\frac{\pi}{2}\big)\big\}$. Then through a~simple symmetry transformation we obtain solutions in the left half plane $S_2 := \big\{w\colon \arg(w) \in \big(\frac{\pi}{2},\frac{3\pi}{2}\big)\big\}$. We start with the formal expansions of the solutions.

Assume that $d$ is a ray of the form ${\rm e}^{{\rm i}\phi} \mathbb{R}^+$ with $\phi \in \big({-}\frac{\pi}{2},\frac{\pi}{2}\big)$. We have the following results on transseries solutions, formal expansions in powers of~$1/w$ and ${\rm e}^{-w}$ (see~\cite{OCInv}), of~\eqref{eqh} valid on~$d$ and, moreover, in the sector~$S_1$:

\begin{Proposition} \label{NFformal}
Assume that $d$ is a ray of the form ${\rm e}^{{\rm i}\phi} \mathbb{R}^+$ with $\phi\in\big({-}\frac{\pi}{2},\frac{\pi}{2}\big)$. Then
\begin{enumerate}[$(i)$]\itemsep=0pt
\item the one-parameter family of transseries solutions of \eqref{eqh} satisfying $h(w) \to 0$ as $|w| \to \infty$ on d are
\begin{gather}\label{transseriesh}
\tilde{h}(w) = \tilde{h}_{0}(w) +\sum_{k=1}^{\infty} C^k {\rm e}^{-kw}w^{-\beta_1k}\tilde{s}_{k}(w),
\end{gather}
where for each $k\geq 1$
\begin{gather*}
\tilde{s}_{k}(w) = \sum_{j=0}^{\infty}\frac{{s}_{k,j}}{w^j}
\end{gather*}
is a formal power series in $w^{-1}$.
\item The formal power series in $w^{-1}$
\begin{gather*}
\tilde{h}_{0}(w)= \sum_{j=2}^{\infty}\frac{{h}_{0,j}}{w^j}
\end{gather*}
is the unique formal power series solution of \eqref{eqh}.
\end{enumerate}
\end{Proposition}
The results in \cite{OCInv} provide us with the relation between these transseries solutions and actual solutions.

In the following, we denote by $\mathcal{L}_{\phi}$ the Laplace transform
\begin{gather*}
f\longmapsto \int_{0}^{\infty {\rm e}^{{\rm i}\phi}} f(p){\rm e}^{-xp}{\rm d}p,
\end{gather*}
where $\phi \in \mathbb{R}$. See also \cite[p.~8]{OCIMRN} for the notation.

\begin{Theorem}\label{NFactual} Let $d$, $\tilde{h}_{0}(w)$ and $\tilde{s}_{k}(w)$ be as in Proposition~{\rm \ref{NFformal}}. Let $h(w)$ be a solution to~\eqref{eqh} on~$d$ for~$|w|$ large enough satisfying
\begin{gather*}
h(w) \to 0,\qquad w\in d, \qquad |w| \to \infty.
\end{gather*}
Then
\begin{enumerate}[$(i)$]\itemsep=0pt
\item There is a unique pair of constants $(C_+,C_-)$ associated with $h(w)$, and $h(w)$ has the following representations
\begin{gather}\label{truesolnup}
h(w) = \mathcal{L}_{\phi} {H}_0(w) + \sum_{k=1}^{\infty}C_{+}^k {\rm e}^{-kw} w^{k M_1} \mathcal{L}_{\phi}{H}_k(w),
\qquad -\phi = \arg(w) \in \left(0,\frac{\pi}{2}\right),\\
\label{truesolndown}
h(w) = \mathcal{L}_{\phi} {H}_0(w) + \sum_{k=1}^{\infty}C_{-}^k {\rm e}^{-kw} w^{k M_1} \mathcal{L}_{\phi}{H}_k(w), \qquad -\phi = \arg(w) \in \left(-\frac{\pi}{2},0\right),
\end{gather}
where
\begin{gather*}
M_1 = \lfloor \operatorname{Re}(-\beta_1)\rfloor+1,\qquad H_0 = \mathcal{B} \tilde{h}_{0},\\
H_k = \mathcal{B} \tilde{h}_{k} = \mathcal{B}\big(w^{-k\beta_1-kM_1}\tilde{s}_{k}\big), \qquad k=1,2,\dots,
\end{gather*}
where each $H_k$ is analytic on the Riemann surface of $\mathbb{C}\backslash\left(\mathbb{Z}^+\cup \mathbb{Z}^-\right)$, and the branch cut for each ${H}_k$, $k\geq 1$, is chosen to be $(-\infty,0]$.

\item There exists $\epsilon_0>0$ such that for each $0<\epsilon \leq \epsilon_0$ there exist $\delta_\epsilon>0$, $R_\epsilon>0$ such that $h(w)$ can be analytically continued to $($at least$)$ the following region
\begin{gather}\label{Sanep}
S_{{\rm an},\epsilon}(h(w)) = S^+_{\epsilon}\cup S^-_{\epsilon},
\end{gather}
where
\begin{gather}
S^{-}_{\epsilon} = \left\{w\colon |w|>R_\epsilon , \, \arg(w) \in \left[-\frac{\pi}{2}-\epsilon,\frac{\pi}{2}-\epsilon\right] \; \mathrm{and}\; |C_{-} {\rm e}^{-w} w^{-\beta_1}| < \delta_\epsilon^{-1} \right\},\nonumber\\
S^{+}_{\epsilon} = \left\{w\colon |w|>R_\epsilon , \, \arg(w) \in \left[-\frac{\pi}{2}+\epsilon,\frac{\pi}{2}+\epsilon \right] \; \mathrm{and}\; |C_{+} {\rm e}^{-w} w^{-\beta_1}| < \delta_\epsilon^{-1}\right\}.\label{Spm}
\end{gather}
Consequently, $h(w)$ is analytic $($at least$)$ in
\begin{gather}\label{San}
S_{\rm an}(h) = \bigcup_{0<\epsilon \leq \epsilon_0}\big(S^{-}_{\epsilon}\cup S^{+}_{\epsilon}\big).
\end{gather}
\item $h(w) \sim h_0(w)$ in $S_1$.
\end{enumerate}
\end{Theorem}

\begin{Note}\label{NFshape}\quad
\begin{enumerate}[(i)]\itemsep=0pt
\item It is straightforward to check that if $\operatorname{Re}(\beta_1) >0$, $S_{\rm an}$ contains all but a compact subset of~${\rm i}\mathbb{R}$. In other words there exists $R_0>0$ such that $h(w)$ is analytic in the closure of~${S}_1 \backslash {\mathbb{D}}_{R_0}$, where $S_1$ is the open right half plane and $\mathbb{D}_{R_0} = \{|w| < R_0\}$ is the open disk centered at origin with radius~$R_0$.
\item On the other hand if $\operatorname{Re}(\beta_1) <0$, $S^c_{\rm an}$ contains all but a compact subset of ${\rm i} \mathbb{R}$. We point out that in particular, the solution is not analytic in ${S}_1 \backslash {\mathbb{D}}_{R_0}$ for any $R_0>0$, contrary to the claim in~\cite{Lin}. Singularities of the tronqu{\'e}e solutions exist for large $w$ in ${S}_1$ as seen in Theorem~\ref{sing},~\ref{NFtri}, \ref{P31sing}, \ref{P31tri}, \ref{P42sing} and \ref{P42tri}.
\end{enumerate}
\end{Note}

\begin{Theorem}[asymptotic position of singularities]\label{sing} Let $h$, $C_{+}$ and $C_{-}$ be as in Theorem~{\rm \ref{NFactual}}.
\begin{enumerate}[$(i)$]\itemsep=0pt
\item Assume $C_{+} \neq 0$. Denote
\begin{gather}\label{xip}
\xi_+(w) = C_{+} w^{-\beta_1} {\rm e}^{-w}.
\end{gather}
Then
\begin{gather*}
h(w) \sim \sum_{m=0}^{\infty} \frac{F_m(\xi_+(w))}{w^m},\qquad |w| \to \infty, \qquad w\in \mathcal{D}^+_w,
\end{gather*}
where for each $m\geq 0$, $F_m$ is analytic at $\xi = 0$ and
\begin{gather}\label{Dwp}
\mathcal{D}^+_w =\left\{|w| > R\colon \arg w \in \left(-\frac{\pi}{2}+\delta,\frac{\pi}{2}+\delta \right), \,\operatorname{dist}\left(\xi_+(w),\Xi \right) >\epsilon ,\, |\xi(w)|<\epsilon^{-1}\right\}\!\!\!
\end{gather}
for any $\delta,\epsilon>0$ small enough and $R$ large enough, and where $\Xi$ is the set of singularities of $F_0(\xi)$. $F_0(\xi)$ satisfies
\begin{gather}\label{conF0}
F_0(0) =0,\qquad F'_0(0)=1.
\end{gather}

\item Assume $C_{+} \neq 0$, and $\xi_s \in \Xi$ is a singularity of $F_0$. Then the singular points of $h$, $w^{+}_n$, near the boundary $\{w\colon \arg(w) = {\pi}/{2} \}$ of the sector of analyticity are given asymptotically by
\begin{gather}\label{wpn}
w^{+}_{n} = 2n\pi{\rm i} -\beta_1 \ln(2n\pi{\rm i})+\ln(C_+)-\ln(\xi_s)+o(1)
\end{gather}
as $n \to \infty$.

\item Assume $C_{-} \neq 0$. Denote
\begin{gather}\label{xim}
\xi_-(w) = C_{-} w^{-\beta_1} {\rm e}^{-w}.
\end{gather}
Then
\begin{gather*}
h(w) \sim \sum_{m=0}^{\infty} \frac{F_m(\xi_-(w))}{w^m},\qquad |w| \to \infty, \qquad w\in \mathcal{D}^-_w,
\end{gather*}
where
\begin{gather}\label{Dwm}
\mathcal{D}^-_w =\left\{|w| > R\colon \arg w \in \left(-\frac{\pi}{2}-\delta,\frac{\pi}{2}-\delta \right), \,\operatorname{dist} (\xi_-(w),\Xi ) >\epsilon , \, |\xi(w)|<\epsilon^{-1}\right\}\!\!\!
\end{gather}
for any $\delta,\epsilon>0$ small enough and $R$ large enough, and where $F_m$, $m\geq 0$, and $\Xi$ are as described in $(i)$.

\item Assume $C_{-} \neq 0$, and $\xi_s \in \Xi$ is a singularity of $F_0$. Then the singular points of $h$, $w^{-}_n$, near the boundary $ \{w\colon \arg(w) = {-\pi}/{2} \}$ of the sector of analyticity are given asymptotically by
\begin{gather*}
w^{-}_{n} = -2n\pi{\rm i} -\beta_1 \ln(-2n\pi{\rm i})+\ln(C_-)-\ln(\xi_s)+o(1)
\end{gather*}
as $n \to \infty$.
\end{enumerate}
\end{Theorem}

The expression of $F_0$ (see \eqref{P31F0}, \eqref{P32F0}, \eqref{P41F0} and \eqref{P42F0}) is obtained explicitly in each case where asymptotic position of singularities is presented.

\subsection{Tritronqu{\'e}e solutions of (\ref{eqh})}\label{NFtrit}

The information on formal and actual tronqu\'{e}e solutions of~\eqref{eqh} in the left half plane $S_2 := \big\{w\colon \arg(w) \in \big(\frac{\pi}{2},\frac{3\pi}{2}\big)\big\}$ is obtained by means of a simple transformation
\begin{gather*}
h(w)=\hat{h}(-w),\qquad \tilde{w} = -w.
\end{gather*}
\eqref{eqh} is rewritten as
\begin{gather} \label{eqhhat}
\hat{h}''(\tilde{w})-\hat{h}(\tilde{w})+\frac{1}{\tilde{w}}\big[(\beta_1-\beta_2) \hat{h}(\tilde{w})+(\beta_1+\beta_2) \hat{h}'(\tilde{w})\big] = g\big({-}\tilde{w},\hat{h},-\hat{h}'\big),
\end{gather}
which is of the form \eqref{eqh} with $\beta_1$ and $\beta_2$ exchanged, and thus all results in Proposition~\ref{NFformal}, Theorems~\ref{NFactual} and~\ref{sing} apply. Without repeating all of the results, we introduce some notations needed for describing the tritronqu\'{e}e solutions of~\eqref{eqh}.

The small transseries solutions of \eqref{eqhhat} in the right half $\tilde{w}$-plane is
\begin{gather*}
\tilde{h}_l(\tilde{w}) = \tilde{h}_{0}(-\tilde{w}) +\sum_{k=1}^{\infty} C^k {\rm e}^{-k\tilde{w}}\tilde{w}^{-\beta_2 k}\tilde{t}_{k}(\tilde{w}),
\end{gather*}
where for each $k\geq 1$, $\tilde{t}_{k}(\tilde{w})$ is a formal power series in $\tilde{w}^{-1}$.

Assume that $\hat{h}(\tilde{w})$ is an actual solution to \eqref{eqhhat} on $d = {\rm e}^{{\rm i}\theta} \mathbb{R}^+$ with $\cos \theta >0$, such that $\hat{h}(\tilde{w})=o(1)$ as $|\tilde{w}| \to \infty$. Then there exists a unique pair of constants $\big(\hat{C}_+,\hat{C}_-\big)$ such that
\begin{gather*}
\hat{h}(\tilde{w}) =
\begin{cases}
\displaystyle \mathcal{L}_{\phi} {\hat{H}}_0(\tilde{w}) + \sum_{k=1}^{\infty}\hat{C}_{+}^k {\rm e}^{-k\tilde{w}} \tilde{w}^{k M_2} \mathcal{L}_{\phi}{\hat{H}}_k(\tilde{w}), & \displaystyle -\phi = \arg(\tilde{w}) \in \left(0,\frac{\pi}{2}\right),\\
\displaystyle \mathcal{L}_{\phi} {\hat{H}}_0(\tilde{w}) + \sum_{k=1}^{\infty}\hat{C}_{-}^k {\rm e}^{-k\tilde{w}} \tilde{w}^{k M_2} \mathcal{L}_{\phi}{\hat{H}}_k(\tilde{w}), & \displaystyle -\phi = \arg(\tilde{w}) \in \left(-\frac{\pi}{2},0\right),
\end{cases}
\end{gather*}
where
\begin{gather}
M_2 = \lfloor \operatorname{Re}(-\beta_2)\rfloor+1,\!\qquad \hat{H}_0(p) = -H_0(-p),\!\qquad
\hat{H}_k = \mathcal{B}\big(w^{-k\beta_2-kM_2}\tilde{t}_{k}\big),\!\qquad k\geq1,\!\!\!\label{hathcdn}
\end{gather}
where each $\hat{H}_k$ is analytic in the Riemann surface of $\mathbb{C}\backslash (\mathbb{Z}^+\cup \mathbb{Z}^-)$, and the branch cut for each~${\hat{H}}_k$, $k\geq 1$, is chosen to be~$(-\infty,0]$. Note that the second equation in~\eqref{hathcdn} holds because the power series solution of~\eqref{eqhhat} must be~$\tilde{h}_0(-\tilde{w})$. By the definition of the Borel transform (see Appendix~\ref{Appx}) we have $\hat{H}_0(p) = -H_0(-p)$.

By Theorem~\ref{NFactual}(ii), $\hat{h}$ is analytic at least on
\begin{gather*}
\hat{S}_{\rm an}(h):=-{S}_{\rm an}\big(\hat{h}\big),
\end{gather*}
where $S_{\rm an}(\hat{h})$ is given by \eqref{Sanep}--\eqref{San} with $\beta_1$ replaced by $\beta_2$. Denote $\hat{\xi}_{\pm} = \hat{C}_{\pm} \tilde{w}^{-\beta_2} {\rm e}^{-\tilde{w}}$ as in~\eqref{xip} and~\eqref{xim}. By Theorem~\ref{sing} if $\hat{C}_+ \neq 0$ then
\begin{gather*}
\hat{h}(\tilde{w}) \sim \sum_{m=0}^{\infty} \frac{\hat{F}_m\big(\hat{\xi}_+(\tilde{w})\big)}{\tilde{w}^m}, \qquad |\tilde{w}| \to \infty, \qquad \tilde{w}\in \mathcal{D}^+_{\tilde{w}},
 \end{gather*}
where $\hat{F}_m$ are analytic at $\xi = 0$. If $\hat{C}_- \neq 0$ then
\begin{gather*}
\hat{h}(\tilde{w}) \sim \sum_{m=0}^{\infty} \frac{\hat{F}_m(\hat{\xi}_-(\tilde{w}))}{\tilde{w}^m}, \qquad |\tilde{w}| \to \infty, \qquad \tilde{w}\in \mathcal{D}^-_{\tilde{w}},
\end{gather*}
where $D^{\pm}_{\tilde{w}}$ are defined by~\eqref{Dwp} and~\eqref{Dwm} with~$\xi_{\pm}$ replaced by $\hat{\xi}_{\pm}$ respectively and~$\Xi$ replaced by~$\hat{\Xi}$ which is defined to be the set of singularities of $\hat{F}_0$.

Tritronqu\'{e}e solutions are special cases of tronqu\'{e}e solutions with $C_{+}=0$ or $C_{-}=0$. Denote
\begin{alignat*}{3}
& h^{+}(w) = \mathcal{L}_{\phi} H_0 (w), \qquad && -\phi = \arg(w) \in (0,\pi),&\\
& h^{-}(w) = \mathcal{L}_{\phi}H_0(w), \qquad && -\phi = \arg(w) \in (-\pi,0),&\\
&\hat{h}^{+}(w)= \mathcal{L}_{\phi} \hat{H}_0 (w), \qquad && -\phi = \arg(w) \in (0,\pi),&\\
&\hat{h}^{-}(w)= \mathcal{L}_{\phi} \hat{H}_0(w), \qquad && -\phi = \arg(w) \in (-\pi,0).&
\end{alignat*}

\begin{Corollary} \label{NFtri}Assume $\phi \in \big(0,\frac{\pi}{2}\big)$. Let ${C}^{t}_{1}$, ${C}^{t}_{2}$, ${C}^{t}_{3}$, ${C}^{t}_{4}$ be the constants in the transseries of~$h^{\pm}$ and~$\hat{h}^{\pm}$, namely,
\begin{gather}
h^{+}(w) = \mathcal{L}_{-\phi} H_0 (w) = \mathcal{L}_{\phi} H_0 (w) + \sum_{k=1}^{\infty}\big({C}^{t}_{1}\big)^k {\rm e}^{-kw} w^{M_1 k} \mathcal{L}_{\phi} H_k (w),\nonumber \\
h^{-}(w) = \mathcal{L}_{\phi} H_0 (w) = \mathcal{L}_{-\phi} H_0 (w) + \sum_{k=1}^{\infty}\big({C}^{t}_{2}\big)^k {\rm e}^{-kw} w^{M_1 k} \mathcal{L}_{-\phi} H_k (w),\nonumber \\
\hat{h}^{+}(w)= \mathcal{L}_{-\phi} \hat{H}_0 (w) = \mathcal{L}_{\phi} \hat{H}_0 (w) + \sum_{k=1}^{\infty}\big({C}^{t}_{3}\big)^k {\rm e}^{-kw} w^{M_2 k} \mathcal{L}_{\phi} \hat{H}_k (w),\nonumber \\
\hat{h}^{-}(w)= \mathcal{L}_{\phi} \hat{H}_0 (w) = \mathcal{L}_{-\phi} \hat{H}_0 (w) + \sum_{k=1}^{\infty}\big({C}^{t}_{4}\big)^k {\rm e}^{-kw} w^{M_2 k} \mathcal{L}_{-\phi} \hat{H}_k (w).\label{htritrep}
\end{gather}

\begin{enumerate}[$(i)$]\itemsep=0pt
\item We have
\begin{gather*}
h^+(w) = \hat{h}^-(-w),\qquad h^-(w) = \hat{h}^+(-w).
\end{gather*}
A consequence of Theorem~{\rm \ref{NFactual}}$(ii)$ is that for any $\delta>0$ there exists $R>0$ such that $h^+$ is analytic in the sector
\begin{gather*}
T^+_{\delta,R}:=\left\{w\colon |w| > R,\, \arg (w) \in \left[-\frac{\pi}{2}+\delta,\frac{3\pi}{2}-\delta\right] \right\},
\end{gather*}
and $h^-$ is analytic in the sector
\begin{gather*}
T^-_{\delta,R}:=\left\{w\colon |w| > R,\, \arg (w) \in \left[-\frac{3\pi}{2}+\delta,\frac{\pi}{2}-\delta\right] \right\}.
\end{gather*}

\item Assume $\xi_s \in \Xi$ is a singularity of $F_0$ $($see Theorem~{\rm \ref{sing}}$(ii))$ and $\hat{\xi}_s \in \hat{\Xi}$ is a singularity of~$\hat{F}_0$. Then the singular points of $h^+$, $w^-_{1,n}$ near the boundary $\big\{w\colon \arg w =- \frac{\pi}{2}\big\}$ and $w^+_{1,n}$ near the boundary $\big\{w\colon \arg w = \frac{3\pi}{2}\big\}$, are given asymptotically by
\begin{gather}
w^-_{1,n} = -2 n \pi {\rm i} -\beta_1 \ln(-2n\pi {\rm i})+\ln \big({C}^{t}_{1}\big) - \ln(\xi_s)+o(1), \nonumber\\
w^+_{1,n} = -2 n \pi {\rm i} +\beta_2 \ln(2n\pi {\rm i})-\ln \big({C}^{t}_{4}\big) +\ln\big(\hat{\xi}_s\big)+o(1),\label{poledistrn}
\end{gather}
as $n \to \infty$. The singular points of $h^-$, $w^-_{2,n}$ near the boundary $\big\{w\colon \arg w =- \frac{3\pi}{2}\big\}$ and~$w^+_{2,n}$ near the boundary $\big\{w\colon \arg w = \frac{\pi}{2}\big\}$, are given asymptotically by
\begin{gather*}
w^-_{2,n} = 2 n \pi {\rm i} +\beta_2 \ln(-2n\pi {\rm i})-\ln \big({C}^{t}_{3}\big) + \ln\big(\hat{\xi}_s\big)+o(1), \\
w^+_{2,n} = 2 n \pi {\rm i} -\beta_1 \ln(2n\pi {\rm i})+\ln \big({C}^{t}_{2}\big) - \ln({\xi}_s)+o(1).
\end{gather*}
\end{enumerate}
\end{Corollary}

\section[Normalizations and Tronqu{\'e}e solutions of $\rm P_{III}$ and $\rm P_{IV}$]{Normalizations and Tronqu{\'e}e solutions of $\boldsymbol{\rm P_{III}}$ and $\boldsymbol{\rm P_{IV}}$}\label{BS3}
\subsection[Tronqu{\'e}e solutions of ${\rm P}^{(i)}_{\rm III}$]{Tronqu{\'e}e solutions of $\boldsymbol{{\rm P}^{(i)}_{\rm III}}$}

If $y(x)$ is a solution of \eqref{P31} which is asymptotic to a formal power series on a ray~$d$ which is not an antistokes line (lines on which $\arg{w} = \pm \frac{\pi}{2}$ where $w$ is the independent variable in the normalized equation), then by dominant balance we have
\begin{gather*}
y(x) \sim l(x), \qquad |x|\to \infty, \qquad x\in d,
\end{gather*}
where
\begin{gather*}
l(x) = A - \left(\frac{\alpha+A^2 \beta}{4} \right)\frac{1}{x}
\end{gather*}
for some $A$ satisfying $A^4=1$. Fix some $A$ satisfying $A^4 = 1$ and make the change of variables
\begin{gather}\label{P31cov}
w=2Ax, \qquad y(x)=h(w)+l\left(\frac{w}{2A}\right).
\end{gather}
Then the equation \eqref{P31} is transformed into an equation for $h$ of the form \eqref{eqh} with
\begin{gather*}
\beta_1 = \frac{1}{2}+\frac{\alpha}{4}-\frac{A^2 \beta}{4}, \qquad \beta_2 = \frac{1}{2}-\frac{\alpha}{4} +\frac{A^2 \beta}{4}.
\end{gather*}
Results in Section~\ref{NFinfo} apply. Let the notations be the same as in Section~\ref{NFinfo}.
\begin{Theorem}\quad
\begin{enumerate}[$(i)$]\itemsep=0pt
\item There is a unique formal power series solution
\begin{gather*}
\tilde{y}_0(x) = \sum_{k=0}^{\infty} \frac{y_{0,k}}{x^k}
\end{gather*}
to \eqref{P31}, where
\begin{gather*}
y_{0,0}=A,\qquad y_{0,1}= -\frac{\alpha+A^2\beta}{4}.
\end{gather*}

\item There is a one-parameter family $\mathcal{F}_{A,1}$ of tronqu{\'e}e solutions of~\eqref{P31} in $A^{-1}S_1$ with representations
\begin{gather}\label{P31FA1}
y(x) = \begin{cases}
\displaystyle l(x)+ h^{+}(2Ax) + \sum_{k=1}^{\infty}C_{+}^k {\rm e}^{-2Akx} (2Ax)^{kM_1} \mathcal{L}_{\phi}{H}_k(2Ax), & \displaystyle -\phi \in \left(0,\frac{\pi}{2}\right),\\
\displaystyle l(x)+ h^{-}(2Ax) + \sum_{k=1}^{\infty}C_{-}^k {\rm e}^{-2Akx} (2Ax)^{kM_1} \mathcal{L}_{\phi}{H}_k(2Ax), & \displaystyle -\phi \in \left(-\frac{\pi}{2},0\right).
\end{cases}\hspace*{-10mm}
\end{gather}

\item There is a one-parameter family $\mathcal{F}_{A,2}$ of tronqu{\'e}e solutions of \eqref{P31} in $A^{-1}S_2$ with representations
\begin{gather}\label{P31FA2}
y(x) =
\begin{cases}
\displaystyle l(x)+ h^{-}(2Ax) + \sum_{k=1}^{\infty}\hat{C}_{+}^k {\rm e}^{2Akx} (-2Ax)^{kM_2} \mathcal{L}_{\phi}{\hat{H}}_k(-2Ax), & \displaystyle -\phi \in \left(0,\frac{\pi}{2}\right),\\
\displaystyle l(x)+ h^{+}(2Ax) + \sum_{k=1}^{\infty}\hat{C}_{-}^k {\rm e}^{2Akx} (-2Ax)^{kM_2} \mathcal{L}_{\phi}{\hat{H}}_k(-2Ax), & \displaystyle -\phi \in \left(-\frac{\pi}{2},0\right).
\end{cases}\hspace*{-10mm}
\end{gather}

\item For each tronqu\'{e}e solution in $(ii)$ or $(iii)$ we have
\begin{gather*}
y(x) \sim \tilde{y}_0(x),\qquad x\in d = A^{-1}{\rm e}^{{\rm i}\theta} \mathbb{R}^+, \qquad |x|\to\infty,
\end{gather*}
and the solution is analytic at least in $(2A)^{-1}{S}_{\rm an}$ if $\cos\theta>0$, in $(2A)^{-1}\hat{S}_{\rm an}$ if $\cos\theta <0$. $S_{\rm an}$ and $\hat{S}_{\rm an}$ are as defined in Theorem~{\rm \ref{NFactual}} and Section~{\rm \ref{NFtrit}}.
\end{enumerate}
\end{Theorem}

From Theorem~\ref{sing} we obtain information about the singularities of $y$. Assume that $y$ is a~tronqu{\'e}e solution with representation~\eqref{P31FA1} or~\eqref{P31FA2}. Let $\xi_+ = C_+ {\rm e}^{-w} w^{-\beta_1}$, $\xi_- = C_- {\rm e}^{-w} w^{-\beta_2}$, $F_m$ and $\hat{F}_m$ be as in Section~\ref{NFinfo}. Then the equation satisfied by $F_0$ is
\begin{gather}\label{P31eqF0}
{\xi}^{2}{\frac {{\rm d}^{2}}{{\rm d}{\xi}^{2}}}F_0 (\xi) +\xi {\frac {{\rm d}}{{\rm d}\xi}}F_0 (\xi) -{\frac {{\xi}^{2} \big( {\frac {{\rm d}}{{\rm d} \xi}}F_0 (\xi) \big) ^{2}}{A+F_0 (\xi) }}- {\frac {( A+F_0 (\xi)) ^{3}}{4{A}^{2}}}+ {\frac {1}{4{A}^{2} ( A+F_0 (\xi)) }} =0.
\end{gather}
The equation of $\hat{F}_0$ is the same as \eqref{P31eqF0}. The solution satisfying \eqref{conF0} is
\begin{gather}\label{P31F0}
F_0(\xi) = \frac{2A \xi}{2A-\xi}.
\end{gather}
\begin{Theorem}\label{P31sing}\quad
\begin{enumerate}[$(i)$]\itemsep=0pt
\item Assume $y(x) \in \mathcal{F}_{A,1}$ is given by the representation~\eqref{P31FA1}. If $C_{+} \neq 0$, then the singular points of $y$, $x^{+}_n$, near the boundary $\{x\colon \arg(2Ax) = {\pi}/{2} \}$ of the sector of analyticity are given asymptotically by
\begin{gather*} 
(2A)x^+_n = 2n\pi {\rm i} -\beta_1 \ln(2n\pi {\rm i})+\ln(C_+)-\ln(2A)+o(1), \qquad n \to \infty.
\end{gather*}

If $C_{-} \neq 0$, then the singular points of $y$, $x^{-}_n$, near the boundary $\{x\colon \arg(2Ax) = {-\pi}/{2}\}$ of the sector of analyticity are given asymptotically by
\begin{gather*}
(2A)x^-_n = -2n\pi {\rm i} -\beta_1 \ln(-2n\pi {\rm i})+\ln(C_{-})-\ln(2A)+o(1), \qquad n \to \infty.
\end{gather*}

\item Assume $y(x) \in \mathcal{F}_{A,2}$ is given by the representation \eqref{P31FA2}. If $\hat{C}_{+} \neq 0$, then the singular points of $y$, $\tilde{x}^{+}_n$, near the boundary $\{\tilde{x}\colon \arg(-2A\tilde{x}) = {\pi}/{2}\}$ of the sector of analyticity are given asymptotically by
\begin{gather*}
(-2A)\tilde{x}^+_n = 2n\pi {\rm i} -\beta_2 \ln(2n\pi {\rm i})+\ln\big(\hat{C}_+\big)-\ln(2A)+o(1), \qquad n \to \infty.
\end{gather*}

If $\hat{C}_{-} \neq 0$, then the singular points of $y$, $\tilde{x}^{-}_n$, near the boundary $\{\tilde{x}\colon \arg(-2A\tilde{x}) = -{\pi}/{2} \}$ of the sector of analyticity are given asymptotically by
\begin{gather*}
(-2A)\tilde{x}^-_n = -2n\pi {\rm i} -\beta_2 \ln(-2n\pi {\rm i})+\ln\big(\hat{C}_{-}\big)-\ln(2A)+o(1), \qquad n \to \infty.
\end{gather*}
\end{enumerate}
\end{Theorem}

From Theorem~\ref{sing} we obtain the following results about tritronqu{\'e}e solutions of \eqref{P31}:
\begin{Theorem}\label{P31tri}
\eqref{P31} has two tritronqu{\'e}e solutions $y^+(x)$ and $y^-(x)$ given by
\begin{gather*}
y^+(x) = l(x) +h^{+}(2Ax),\qquad y^-(x)= l(x)+h^{-}(2Ax).
\end{gather*}
Let $C^t_{j}$, $1\leq j\leq 4$ be as in \eqref{htritrep}. Then
\begin{enumerate}[$(i)$]\itemsep=0pt
\item $\mathcal{F}_{A,1} \cap \mathcal{F}_{A,2} = \{ y^+,y^-\}$.

\item For each $\delta>0$ there exists $R>0$ such that $y^+(x)$ is analytic in $A^{-1}T^+_{\delta,R}$, and $y^+$ is asymptotic to $y_0(x)$ in the sector
\begin{gather*}
\bigcup_{-\frac{\pi}{2}< \phi < \frac{3\pi}{2} } \big(A^{-1}{\rm e}^{{\rm i}\phi} \mathbb{R}^+\big) .
\end{gather*}

The singular points of $y^+(x)$, $x^{\pm}_{1,n}$, near the boundary of the sector of analyticity are given asymptotically by
\begin{gather*}
(2A)x^-_{1,n} = -2n\pi {\rm i} -\beta_1 \ln(-2n\pi {\rm i})+\ln\big(C^{t}_{1}\big)-\ln(2A)+o(1), \qquad n \to \infty,\\
(2A){x}^+_{1,n} = -2n\pi {\rm i} +\beta_2 \ln(2n\pi {\rm i})-\ln\big({C}^{t}_4\big)+\ln(2A)+o(1), \qquad n \to \infty.
\end{gather*}

\item For each $\delta>0$ there exists $R>0$ such that $y^-(x)$ is analytic in $A^{-1}T^-_{\delta,R}$, and $y^-$ is asymptotic to $y_0(x)$ in the sector
\begin{gather*}
\bigcup_{-\frac{3\pi}{2}< \phi < \frac{\pi}{2} } \big(A^{-1}{\rm e}^{{\rm i}\phi} \mathbb{R}^+\big).
\end{gather*}

The singular points of $y^-(x)$, $x^{\pm}_{2,n}$, near the boundary of the sector of analyticity are given asymptotically by
\begin{gather*}
(2A){x}^-_{2,n} = 2n\pi {\rm i} +\beta_2 \ln(-2n\pi {\rm i})-\ln\big({C}^{t}_3\big)+\ln(2A)+o(1), \qquad n \to \infty,\\
(2A)x^+_{2,n} = 2n\pi {\rm i} -\beta_1 \ln(2n\pi {\rm i})+\ln\big(C^{t}_{2}\big)-\ln(2A)+o(1), \qquad n \to \infty.
\end{gather*}
\end{enumerate}
\end{Theorem}

\subsection[Tronqu{\'e}e solutions of ${\rm P}^{(ii)}_{\rm III}$]{Tronqu{\'e}e solutions of $\boldsymbol{{\rm P}^{(ii)}_{\rm III}}$}

If $y(x)$ is a solution of \eqref{P32} which is asymptotic to a formal power series on a ray $d$ which is not an antistokes line, then by dominant balance we have
\begin{gather*}
y(x) \sim l(x), \qquad |x|\to \infty, \qquad x\in d,
\end{gather*}
where
\begin{gather}\label{P32lx}
l(x) = A x^{1/3} -\frac{\beta}{3A x^{1/3}}
\end{gather}
for some $A$ satisfying $A^3=1$. Fix an $A$ satisfying $A^3 = 1$ and make the change of variables
\begin{gather}\label{P32cov}
w=(27A/4)^{1/2}x^{2/3}, \qquad y(x)=x^{1/3}h(w)+l(x).
\end{gather}
Then the equation \eqref{P32} is transformed into an equation for $h$ of the form \eqref{eqh} with
\begin{gather*}
\beta_1 =\beta_2 = \frac{1}{2}.
\end{gather*}
Let the notations be the same as in Section~\ref{NFinfo}. In view of the transformation \eqref{P32cov}, we denote
\begin{gather*}
S^{(0)}_{R}:=\left\{x\colon |x|\geq R, \, \arg(x) \in \left[-\frac{3\pi}{4}-\frac{3\arg{{A}}}{4},\frac{3\pi}{4}-\frac{3\arg{{A}}}{4}\right]\right\}, \\
S^{(1)}_{R}:=\left\{x\colon |x|\geq R, \,\arg(x) \in \left[\frac{3\pi}{4}-\frac{3\arg{{A}}}{4},\frac{9\pi}{4}-\frac{3\arg{{A}}}{4}\right]\right\},\\
S^{(2)}_{R}:=\left\{x\colon |x|\geq R,\, \arg(x) \in \left[\frac{9\pi}{4}-\frac{3\arg{{A}}}{4},\frac{15\pi}{4}-\frac{3\arg{{A}}}{4}\right]\right\},\\
S^{(3)}_{R}:=\left\{x\colon |x|\geq R,\, \arg(x) \in \left[\frac{15\pi}{4}-\frac{3\arg{{A}}}{4},\frac{21\pi}{4}-\frac{3\arg{{A}}}{4}\right]\right\}.
\end{gather*}
We notice that for $j \in\{ 0,2\}$, $S^{(j)}_R$ is mapped under the transformation~\eqref{P32cov} bijectively to the closed sector $\overline{S}_1 \backslash {\mathbb{D}}_{R_0}$ in the $w$-plane, where $R_0 = R^{2/3}$ (see also Note~\ref{NFshape}); for $j \in\{ 1,3\}$, $S^{(j)}_R$ is mapped bijectively to the closed sector $\overline{S}_2 \backslash {\mathbb{D}}_{R_0}$ in the $w$-plane.

\begin{Theorem}\quad
\begin{enumerate}[$(i)$]\itemsep=0pt
\item There is a unique formal power series solution
\begin{gather*}
\tilde{y}_0(x) = x^{1/3} \sum_{k=0}^{\infty} \frac{y_{0,k}}{x^{2k/3}}
\end{gather*}
to \eqref{P32}, where
\begin{gather*}
y_{0,0}=A,\qquad y_{0,1}= -\frac{\beta}{3A}.
\end{gather*}

\item For each $j \in \{0,1,2,3\}$, there is a one-parameter family $\mathcal{F}_{A,j}$ of tronqu{\'e}e solutions of~\eqref{P32} in~$S^{(j)}_R$ where
\begin{gather*}
y(x) = l(x)+ x^{1/3} h(w), \qquad w = K x^{2/3},\qquad K=(27A/4)^{1/2}.
\end{gather*}
If $j$ is even, then $h(w)$ has the representations
\begin{gather}\label{P32FA02}
h(w) = \begin{cases}
\displaystyle h^{+}(w) + \sum_{k=1}^{\infty}C_{+}^k {\rm e}^{-kw} \mathcal{L}_{\phi}{H}_k(w), &\displaystyle -\phi \in \left(0,\frac{\pi}{2}\right],\\
\displaystyle h^{-}(w) + \sum_{k=1}^{\infty}C_{-}^k {\rm e}^{-kw} \mathcal{L}_{\phi}{H}_k(w), & \displaystyle -\phi \in \left[-\frac{\pi}{2},0\right).
\end{cases}
\end{gather}
If $j$ is odd, then $h(w)$ has the representations
\begin{gather}\label{P32FA13}
h(w) = \begin{cases}
\displaystyle h^{-}(w) + \sum_{k=1}^{\infty}\hat{C}_{+}^k {\rm e}^{kw} \mathcal{L}_{\phi}{\hat{H}}_k(-w), & \displaystyle -\phi \in \left(0,\frac{\pi}{2}\right],\\
\displaystyle h^{+}(w) + \sum_{k=1}^{\infty}\hat{C}_{-}^k {\rm e}^{kw} \mathcal{L}_{\phi}{\hat{H}}_k(-w), &\displaystyle -\phi \in \left[-\frac{\pi}{2},0\right).
\end{cases}
\end{gather}

\item Let $y(x)$ be a tronqu{\'e}e solution in $\mathcal{F}_{A,j}$. If $j$ is even, the region of analyticity contains the corresponding branch of $\big(K^{-1} S_{\rm an}(h)\big)^{2/3}$, which contains~$S^{(j)}_R$ for~$R$ large enough. If~$j$ is odd, the region of analyticity contains the corresponding branch of $\big(K^{-1} \hat{S}_{\rm an}(h)\big)^{2/3}$, which contains~$S^{(j)}_R$ for $R$ large enough, and
\begin{gather*}
y(x) \sim \tilde{y}_0(x),\qquad x\in d, \qquad |x|\to\infty,
\end{gather*}
where $d$ is a ray whose infinite part is contained in the interior of~$S^{(j)}_R$.
\end{enumerate}
\end{Theorem}

Assume that $y(x)$ is a tronqu{\'e}e solution to \eqref{P32} and $h$ is defined by~\eqref{P32cov}. Then $h$ has the representation~\eqref{P32FA02} or~\eqref{P32FA13}. From Theorem~\ref{sing} we obtain information about singularities of $h$. Let $\xi_+ = C_+ {\rm e}^{-w} w^{-1/2}$, $\xi_- = C_- {\rm e}^{-w} w^{-1/2}$, $F_m$ and $\hat{F}_m$ be as in Section~\ref{NFinfo}. Then the equation satisfied by $F_0$ and $\hat{F}_0$ is the same
\begin{gather}\label{P32eqF0}
{\xi}^{2}{\frac {{\rm d}^{2}}{{\rm d}{\xi}^{2}}}F_0 (\xi) + \left( {
\frac {{\rm d}}{{\rm d}\xi}}F_0 (\xi) \right) \xi-{\frac { \big( {
\frac {{\rm d}}{{\rm d}\xi}}F_0 (\xi) \big) ^{2}{\xi}^{2}}{{A}+
F_0 (\xi) }}- {\frac {( {A}+F_0 (\xi)) ^{2}}{3{A}}}+ {\frac {1}{3{A} ( {A}+F_0 (\xi) ) }} =0.
\end{gather}
The solution satisfying \eqref{conF0} is
\begin{gather}\label{P32F0}
F_0(\xi) = \frac{36A^2 \xi}{(6A-\xi)^2}.
\end{gather}
\begin{Theorem} \label{P32sing}\quad
\begin{enumerate}[$(i)$]\itemsep=0pt
\item If $j\in\{0,2\}$, then $h$ has representation \eqref{P32FA02} for a unique pair of constants $(C_+,C_-)$. If $C_{+} \neq 0$, then the singular points of $h$, $w^{+}_n$, near the boundary $\{w\colon \arg w = {\pi}/{2}\}$ of the sector of analyticity are given asymptotically by
\begin{gather*}
w^{+}_{n} = 2n\pi {\rm i} -\frac{\ln(2n\pi {\rm i})}{2} +\ln(C_+)-\ln(6A)+o(1), \qquad n \to \infty.
\end{gather*}
If $C_{-} \neq 0$, then the singular points of $h$, $w^{-}_n$, near the boundary $\{w\colon \arg w = -{\pi}/{2}\}$ of the sector of analyticity are given asymptotically by
\begin{gather*}
w^{-}_{n} = -2n\pi {\rm i} -\frac{\ln(-2n\pi {\rm i})}{2}+\ln(C_-)-\ln(6A)+o(1),\qquad n \to \infty.
\end{gather*}

\item If $j\in\{1,3\}$, then $h$ has representation \eqref{P32FA13} for a unique pair of constants $\big(\hat{C}_+,\hat{C}_-\big)$. If $\hat{C}_{+} \neq 0$, then the singular points of $h$, $\tilde{w}^{+}_n$, near the boundary $\{\tilde{w}\colon \arg \tilde{w} = -{\pi}/{2}\}$ of the sector of analyticity are given asymptotically by
\begin{gather*}
\tilde{w}^+_n = -2n\pi {\rm i} +\frac{\ln(2n\pi {\rm i})}{2} -\ln\big(\hat{C}_+\big)+\ln(6A)+o(1), \qquad n \to \infty.
\end{gather*}

If $\hat{C}_{-} \neq 0$, then the singular points of~$h$, $\tilde{w}^{-}_n$, near the boundary $\{\tilde{w}\colon \arg \tilde{w} = {\pi}/{2} \}$ of the sector of analyticity are given asymptotically by
\begin{gather*}
\tilde{w}^{-}_{n} = 2n\pi {\rm i} +\frac{\ln(-2n\pi {\rm i})}{2}-\ln\big(\hat{C}_-\big)+\ln(6A)+o(1), \qquad n \to \infty.
\end{gather*}
\end{enumerate}
\end{Theorem}

\begin{Theorem}\label{P32tri}\quad
\begin{enumerate}[$(i)$]\itemsep=0pt
\item For each $j\in\{0,2\}$ we have a tritronqu{\'e}e solution $y^+_{j}$ analytic in $S^{(j)}_R \bigcup S^{(j+1)}_R$ for $R$ large enough, given by
\begin{gather*}
y^+_{j}(x) = l(x) +h^{+}(w).
\end{gather*}

\item For each $j\in\{1,3\}$ we have a tritronqu{\'e}e solution $y^-_{j}$ analytic in $S^{(j)}_R \bigcup S^{(j+1)}_R$, where $S^{(4)}_R = S^{(0)}_R$ and $R$ is large enough, given by
\begin{gather*}
y^-_{j}(x) = l(x) +h^{-}(w).
\end{gather*}
\end{enumerate}
Let $C^t_{j}$, $1\leq j\leq 4$ be as in \eqref{htritrep}. Then
\begin{enumerate}[(i)]\itemsep=0pt
\item[$(iii)$] The singular points of $h^+(w)$, $w^{\pm}_{1,n}$, near the boundary of the sector of analyticity are given asymptotically by
\begin{gather*}
w^-_{1,n} = -2 n \pi {\rm i} -\frac{\ln(-2n\pi {\rm i})}{2}+\ln \big({C}^{t}_{1}\big) - \ln(6A)+o(1), \\
w^+_{1,n} = -2 n \pi {\rm i} +\frac{ \ln(2n\pi {\rm i})}{2}-\ln \big({C}^{t}_{4}\big) +\ln(6A)+o(1).
\end{gather*}

\item[$(iv)$] The singular points of $h^-(w)$, $w^{\pm}_{2,n}$, near the boundary of the sector of analyticity are given asymptotically by
\begin{gather*}
w^-_{2,n} = 2 n \pi {\rm i} +\frac{ \ln(-2n\pi {\rm i})}{2}-\ln \big({C}^{t}_{3}\big) + \ln(6A)+o(1), \\
w^+_{2,n} = 2 n \pi {\rm i} - \frac{\ln(2n\pi {\rm i})}{2}+\ln \big({C}^{t}_{2}\big) - \ln(6A)+o(1).
\end{gather*}
\end{enumerate}
\end{Theorem}

\subsection[Tronqu{\'e}e solutions of $\rm P_{IV}$]{Tronqu{\'e}e solutions of $\boldsymbol{\rm P_{IV}}$}

By dominant balance we have four possibilities for the leading behavior of $\rm P_{IV}$. We shall study them one by one.
\begin{gather*}
y(x) \sim l(x), \qquad |x|\to \infty, \qquad x\in d.
\end{gather*}

\subsubsection{Case 1}
\begin{gather*}
l(x) = -\frac{2x}{3} + \frac{\alpha}{x}.
\end{gather*}
Make the change of variables
\begin{gather}\label{P41cov}
x=\big(\sqrt{3} {\rm i} w\big)^{1/2}, \qquad y(x)=x h(w)+l(x).
\end{gather}
Then the equation \eqref{P4} is transformed into an equation for $h$ of the form \eqref{eqh} with
\begin{gather*}
\beta_1 = \beta_2 = \frac{1}{2}.
\end{gather*}

Let the notations be the same as in Section~\ref{NFinfo}. In view of the transformation \eqref{P41cov}, we denote
\begin{gather}
S^{(0)}_{R}:=\left\{x\colon |x|\geq R,\, \arg(x) \in \left[0,\frac{\pi}{2}\right]\right\}, \qquad S^{(1)}_{R}:=\left\{x\colon |x|\geq R,\, \arg(x) \in \left[\frac{\pi}{2},\pi \right]\right\},\!\!\!\label{P41SjR}\\
S^{(2)}_{R}:=\left\{x\colon |x|\geq R,\, \arg(x) \in \left[\pi,\frac{3\pi}{2}\right]\right\},\qquad S^{(3)}_{R}:=\left\{x\colon |x|\geq R,\, \arg(x) \in \left[\frac{3\pi}{2},2 \pi\right]\right\}.\nonumber
\end{gather}
We notice that for $j \in\{ 0,2\}$, $S^{(j)}_R$ is mapped under the transformation~\eqref{P32cov} bijectively to the closed sector $\overline{S}_1 \backslash {\mathbb{D}}_{R^2}$ in the $w$-plane, (see also Note~\ref{NFshape}); for $j \in\{ 1,3\}$, $S^{(j)}_R$ is mapped bijectively to the closed sector~$\overline{S}_2 \backslash {\mathbb{D}}_{R^2}$ in the $w$-plane.

\begin{Theorem}\quad 
\begin{enumerate}[$(i)$]\itemsep=0pt
\item There is a formal power series solution of \eqref{P4} of the form
\begin{gather*}
\tilde{y}_0(x) = x \sum_{k=0}^{\infty} \frac{y_{0,k}}{x^{2k}},
\end{gather*}
where
\begin{gather*}
y_{0,0}=-\frac{2}{3},\qquad y_{0,1}= {\alpha}.
\end{gather*}

\item For each $j \in \{0,1,2,3\}$, there is a one-parameter family $\mathcal{F}_{A,j}$ of tronqu{\'e}e solutions of~\eqref{P4} in~$S^{(j)}_R$, where
\begin{gather*}
y(x) = l(x)+ x h(w), \qquad w=\frac{x^2}{\sqrt{3}{\rm i}}.
\end{gather*}
If $j$ is even, then $h(w)$ has the representations
\begin{gather}\label{P41FA02}
h(w) =
\begin{cases}
\displaystyle h^{+}(w) + \sum\limits_{k=1}^{\infty}C_{+}^k {\rm e}^{-kw} \mathcal{L}_{\phi}{H}_k(w), & \displaystyle -\phi \in \left(0,\frac{\pi}{2}\right],\\
\displaystyle h^{-}(w) + \sum\limits_{k=1}^{\infty}C_{-}^k {\rm e}^{-kw} \mathcal{L}_{\phi}{H}_k(w), & \displaystyle -\phi \in \left[-\frac{\pi}{2},0\right).
\end{cases}
\end{gather}
If $j$ is odd, then $h(w)$ has the representations
\begin{gather}\label{P41FA13}
h(w) =
\begin{cases}
\displaystyle h^{-}(w) + \sum_{k=1}^{\infty}\hat{C}_{+}^k {\rm e}^{kw} \mathcal{L}_{\phi}{\hat{H}}_k(-w), & \displaystyle -\phi \in \left(0,\frac{\pi}{2}\right],\\
\displaystyle h^{+}(w) + \sum_{k=1}^{\infty}\hat{C}_{-}^k {\rm e}^{kw} \mathcal{L}_{\phi}{\hat{H}}_k(-w), & \displaystyle -\phi \in \left[-\frac{\pi}{2},0\right).
\end{cases}
\end{gather}

\item Let $y(x)$ be a tronqu{\'e}e solution in $\mathcal{F}_{A,j}$. If $j$ is even, the region of analyticity contains the corresponding branch of $\big(\sqrt{3}{\rm i} S_{\rm an}(h)\big)^{1/2}$, which contains~$S^{(j)}_R$ for~$R$ large enough. If~$j$ is odd, the region of analyticity contains the corresponding branch of $\big(\sqrt{3}{\rm i} \hat{S}_{\rm an}(h)\big)^{1/2}$, which contains~$S^{(j)}_R$ for~$R$ large enough, and
\begin{gather*}
y(x) \sim \tilde{y}_0(x),\qquad x\in d, \qquad |x|\to\infty,
\end{gather*}
where $d$ is a ray whose infinite part is contained in the interior of $S^{(j)}_R$.
\end{enumerate}
\end{Theorem}

Assume that $y(x)$ is a tronqu{\'e}e solution to \eqref{P4} satisfying $y(x) \sim -\frac{2x}{3}$ and $h$ is defined by~\eqref{P41cov}. Then $h$ has representation~\eqref{P41FA02} or~\eqref{P41FA13}. From Theorem~\ref{sing} we obtain information about singularities of $h$. Let $\xi_+ = C_+ {\rm e}^{-w} w^{-1/2}$, $\xi_- = C_- {\rm e}^{-w} w^{-1/2}$, $F_m$ and $\hat{F}_m$ be as in Section~\ref{NFinfo}. Then the equation satisfied by $F_0$ and $\hat{F}_0$ is the same
\begin{gather}
{\xi}^{2}{\frac {{\rm d}^{2}}{{\rm d}{\xi}^{2}}}F_0 (\xi) +\xi {
\frac {{\rm d}}{{\rm d}\xi}}F_0 (\xi) - {\frac {3 {\xi}^{2} \big( {
\frac {{\rm d}}{{\rm d}\xi}}F_0 (\xi) \big) ^{2}}{2(3 F_0 (\xi) -2)}}\nonumber\\
\qquad{} + \frac{( 3 F_0 (\xi) -2) ^{3}}{24}+ \frac{( 3 F_0 (\xi) -2) ^{2}}{3}+\frac{3 F_0 ( \xi)-2}{2} =0.\label{P41eqF0}
\end{gather}
The solution satisfying \eqref{conF0} is
\begin{gather}\label{P41F0}
F_0(\xi) = \frac{4 \xi}{\xi^2+2\xi+4}
\end{gather}
with simple poles at $\xi^{(1)}_s =-1-\sqrt{3}{\rm i}$ and $\xi^{(2)}_s =-1+\sqrt{3}{\rm i}$. Hence the statements in Theorems~\ref{P32sing} and~\ref{P32tri} hold true for function~$h$, with $S^{(j)}_R$ as defined in \eqref{P41SjR} and $6A$ in the formula replaced by $\xi^{(i)}_s$, $i=1$ or $i=2$.

\subsubsection{Case 2}
\begin{gather*}
l(x) = -{2x} - \frac{\alpha}{x}.
\end{gather*}
Make the change of variables
\begin{gather}\label{P42cov}
x=\left( w\right)^{1/2}, \qquad y(x)=x h(w)+l(x).
\end{gather}
Then the equation \eqref{P4} is transformed into an equation for $h$ of the form \eqref{eqh} with
\begin{gather*}
\beta_1 =\alpha+\frac{1}{2}, \qquad \beta_2 = -\alpha+\frac{1}{2}.
\end{gather*}

Let the notations be the same as in Section~\ref{NFinfo}. In view of the transformation~\eqref{P41cov}, we denote
\begin{gather}
S^{(0)} :=\left\{x\colon \arg(x) \in \left(-\frac{\pi}{4},\frac{\pi}{4}\right)\right\}, \qquad S^{(1)}:=\left\{x\colon \arg(x) \in \left(\frac{\pi}{4},\frac{3\pi}{4} \right)\right\},\nonumber\\
S^{(2)} :=\left\{x\colon \arg(x) \in \left(\frac{3\pi}{4},\frac{5\pi}{4}\right)\right\}, \qquad S^{(3)}:=\left\{x\colon \arg(x) \in \left(\frac{5\pi}{4},\frac{7\pi}{4}\right)\right\}.\label{P42SjR}
\end{gather}
For $j \in\{ 0,2\}$, $S^{(j)}$ is mapped under the transformation~\eqref{P32cov} bijectively to the right half $w$-plane~${S}_1$; for $j \in\{ 1,3\}$, $S^{(j)}$ is mapped bijectively to the sector to the left half $w$-plane ${S}_2$.

\begin{Theorem}\quad
\begin{enumerate}[$(i)$]\itemsep=0pt
\item There is a formal power series solution of \eqref{P4} of the form
\begin{gather*}
\tilde{y}_0(x) = x \sum_{k=0}^{\infty} \frac{y_{0,k}}{x^{2k}},
\end{gather*}
where
\begin{gather*}
y_{0,0}=-2,\qquad y_{0,1}= -{\alpha}.
\end{gather*}

\item For each $j \in \{0,1,2,3\}$, there is a one-parameter family $\mathcal{F}_{A,j}$ of tronqu{\'e}e solutions of~\eqref{P4} in $S^{(j)}$, where
\begin{gather*}
y(x) = l(x)+ x h(w), \qquad w=x^2.
\end{gather*}
If $j$ is even, then $h(w)$ has the representations
\begin{gather}\label{P42FA02}
h(w) =
\begin{cases}
\displaystyle h^{+}(w) + \sum_{k=1}^{\infty}C_{+}^k {\rm e}^{-kw} w^{kM_1} \mathcal{L}_{\phi}{H}_k(w), & \displaystyle -\phi \in \left(0,\frac{\pi}{2}\right),\\
\displaystyle h^{-}(w) + \sum_{k=1}^{\infty}C_{-}^k {\rm e}^{-kw} w^{kM_1} \mathcal{L}_{\phi}{H}_k(w), & \displaystyle -\phi \in \left(-\frac{\pi}{2},0\right).
\end{cases}
\end{gather}
If $j$ is odd, then $h(w)$ has the representations
\begin{gather}\label{P42FA13}
h(w) = \begin{cases}
\displaystyle h^{-}(w) + \sum_{k=1}^{\infty}\hat{C}_{+}^k {\rm e}^{kw} (-w)^{kM_2} \mathcal{L}_{\phi}{\hat{H}}_k(-w), &\displaystyle -\phi \in \left(0,\frac{\pi}{2}\right),\\
\displaystyle h^{+}(w) + \sum_{k=1}^{\infty}\hat{C}_{-}^k {\rm e}^{kw} (-w)^{kM_2} \mathcal{L}_{\phi}{\hat{H}}_k(-w), & \displaystyle -\phi \in \left(-\frac{\pi}{2},0\right).
\end{cases}
\end{gather}

\item Let $y(x)$ be a tronqu{\'e}e solution in $\mathcal{F}_{A,j}$. If $j$ is even, then the region of analyticity contains the corresponding branch of $(S_{\rm an}(h))^{1/2}$. If $j$ is odd, then the region of analyticity contains the corresponding branch of $\big(\hat{S}_{\rm an}(h)\big)^{1/2}$, and
\begin{gather*}
y(x) \sim \tilde{y}_0(x),\qquad\ x\in d \subset S^{(j)}, \qquad |x|\to\infty.
\end{gather*}
\end{enumerate}
\end{Theorem}

Assume that $y(x)$ is a tronqu{\'e}e solution to \eqref{P4} satisfying $y(x) \sim -2 x$ and $h$ is defined by~\eqref{P42cov}. Then $h$ has representation~\eqref{P42FA02} or~\eqref{P42FA13}. From Theorem~\ref{sing} we obtain information about singularities of $h$. Let $F_m$ and $\hat{F}_m$ be as in Section~\ref{NFinfo}. Then the equation satisfied by $F_0$ and $\hat{F}_0$ is the same
\begin{gather}
{\xi}^{2}{\frac {{\rm d}^{2}}{{\rm d}{\xi}^{2}}}F_0 (\xi) +\xi {
\frac {{\rm d}}{{\rm d}\xi}}F_0 (\xi) - {\frac {{\xi}^{2} \big( {
\frac {{\rm d}}{{\rm d}\xi}}F_0 (\xi) \big) ^{2}}{2(F_0 ( \xi) -2)}}-\frac{3 ( F_0 (\xi) -2 ) ^{3}}{8}\nonumber\\
\qquad{} -( F_0 (\xi) -2) ^{2}-\frac{ F_0 (\xi)-2}{2}=0.\label{P42eqF0}
\end{gather}
The solution satisfying \eqref{conF0} is
\begin{gather}\label{P42F0}
F_0(\xi) = \frac{2 \xi}{\xi+2}
\end{gather}
with a simple pole at $\xi_s =-2 $.

\begin{Theorem} \label{P42sing}\quad
\begin{enumerate}[$(i)$]\itemsep=0pt
\item If $j\in\{0,2\}$, then $h$ has representation \eqref{P42FA02} for a unique pair of constants $(C_+,C_-)$. If $C_{+} \neq 0$, then the singular points of~$h$, at $w^{+}_n$, near the boundary $\{w\colon \arg w ={\pi}/{2}\big\}$ of the sector of analyticity are given asymptotically by
\begin{gather*}
w^{+}_{n} = 2n\pi {\rm i} -(\alpha+1/2){\ln(2n\pi {\rm i})} +\ln(C_+)-\ln(-2)+o(1), \qquad n \to \infty.
\end{gather*}
If $C_{-} \neq 0$, then the singular points of $h$, $w^{-}_n$, near the boundary $\{w\colon \arg w =-{\pi}/{2}\}$ of the sector of analyticity are given asymptotically by
\begin{gather*}
w^{-}_{n} = -2n\pi {\rm i} - (\alpha+1/2 ){\ln(-2n\pi {\rm i})}+\ln(C_-)-\ln(-2)+o(1),\qquad n \to \infty.
\end{gather*}

\item If $j\in\{1,3\}$, then $h$ has representation~\eqref{P42FA13} for a unique pair of constants $\big(\hat{C}_+,\hat{C}_-\big)$. If $\hat{C}_{+} \neq 0$, then the singular points of~$h$, $\tilde{w}^{+}_n$, near the boundary~$\{\tilde{w}\colon \arg \tilde{w} =-{\pi}/{2}\}$ of the sector of analyticity are given asymptotically by
\begin{gather*}
\tilde{w}^+_n = -2n\pi {\rm i} +(-\alpha+1/2){\ln(2n\pi {\rm i})}-\ln\big(\hat{C}_+\big)+\ln(-2)+o(1), \qquad n \to \infty.
\end{gather*}

If $\hat{C}_{-} \neq 0$, then the singular points of~$h$, $\tilde{w}^{-}_n$, near the boundary $\{\tilde{w}\colon \arg \tilde{w} ={\pi}/{2}\}$ of the sector of analyticity are given asymptotically by
\begin{gather*}
\tilde{w}^{-}_{n} = 2n\pi {\rm i} +(-\alpha+1/2){\ln(-2n\pi {\rm i})}-\ln\big(\hat{C}_-\big)+\ln(-2)+o(1), \qquad n \to \infty.
\end{gather*}
\end{enumerate}
\end{Theorem}

Denote
\begin{gather}
T^{(0)}_{\delta,R}:=\left\{w\colon |w| > R,\, \arg (w) \in \left[-\frac{\pi}{4}+\delta,\frac{3\pi}{4}-\delta\right] \right\},\nonumber\\
T^{(1)}_{\delta,R}:=\left\{w\colon |w| > R,\, \arg (w) \in \left[\frac{\pi}{4}+\delta,\frac{5\pi}{4}-\delta\right] \right\},\nonumber\\
T^{(2)}_{\delta,R}:=\left\{w\colon |w| > R,\, \arg (w) \in \left[\frac{3\pi}{4}+\delta,\frac{7\pi}{4}-\delta\right] \right\},\nonumber\\
T^{(3)}_{\delta,R}:=\left\{w\colon |w| > R,\, \arg (w) \in \left[-\frac{3\pi}{4}+\delta,\frac{\pi}{4}-\delta\right] \right\}.\label{P42Tj}
\end{gather}

\begin{Theorem}\label{P42tri}\quad{}
\begin{enumerate}[$(i)$]\itemsep=0pt
\item Let $j\in\{0,2\}$. For each $\delta>0$ there exists $R$ large enough such that we have a tritronqu{\'e}e solution $y^+_{j}$ analytic in $T^{(j)}_{\delta,R}$ given by
\begin{gather*}
y^+_{j}(x) = -2x+\frac{\alpha}{x} +h^{+}\big(x^2\big).
\end{gather*}

\item Let $j\in\{1,3\}$. For each $\delta>0$ there exists $R$ large enough such that we have a tritronqu{\'e}e solution $y^-_{j}$ analytic in $T^{(j)}_{\delta,R}$ given by
\begin{gather*}
y^-_{j}(x) = -2x+\frac{\alpha}{x} +h^{-}\big(x^2\big).
\end{gather*}
\end{enumerate}
Let $C^t_{j}$, $1\leq j\leq 4$ be as in \eqref{htritrep}. Then
\begin{enumerate}\itemsep=0pt
\item[$(iii)$] The singular points of $h^+(w)$, $w^{\pm}_{1,n}$, near the boundary of the sector of analyticity are given asymptotically by
\begin{gather*}
w^-_{1,n} = -2 n \pi {\rm i} -(\alpha+1/2){\ln(-2n\pi {\rm i})}+\ln ({C}^{t}_{1}) - \ln(-2)+o(1), \\
w^+_{1,n} = -2 n \pi {\rm i} +(-\alpha+1/2){ \ln(2n\pi {\rm i})}-\ln ({C}^{t}_{4}) +\ln(-2)+o(1).
\end{gather*}

\item[$(iv)$] The singular points of $h^-(w)$, $w^{\pm}_{2,n}$, near the boundary of the sector of analyticity are given asymptotically by
\begin{gather*}
w^-_{2,n} = 2 n \pi {\rm i} +(-\alpha+1/2){ \ln(-2n\pi {\rm i})}-\ln ({C}^{t}_{3}) + \ln(-2)+o(1), \\
w^+_{2,n} = 2 n \pi {\rm i} - (\alpha+1/2){\ln(2n\pi {\rm i})}+\ln ({C}^{t}_{2}) - \ln(-2)+o(1).
\end{gather*}
\end{enumerate}
\end{Theorem}

\subsubsection{Case 3}
\begin{gather*}
l(x) = \frac{A}{ x}+ \frac{\alpha A+\beta }{2x^3}, \qquad A^2 = -b/2.
\end{gather*}
Make the change of variables
\begin{gather}\label{P43cov}
x=\left( w\right)^{1/2}, \qquad y(x)=x^{-1} h(w)+l(x).
\end{gather}
Then the equation \eqref{P4} is transformed into an equation for $h$ of the form \eqref{eqh} with
\begin{gather*}
\beta_1 =-\frac{\alpha}{2}+\frac{3A}{2},\qquad \beta_2 = \frac{\alpha}{2}-\frac{3A}{2}.
\end{gather*}

Let the notations be as in Section~\ref{NFinfo}, $S^{(j)}$ be as in \eqref{P42SjR} and $T^{(j)}_{\delta,R}$ be as in~\eqref{P42Tj}.

\begin{Theorem}\quad
\begin{enumerate}[$(i)$]\itemsep=0pt
\item There is a formal power series solution of \eqref{P4} of the form
\begin{gather*}
\tilde{y}_0(x) = \frac{1}{x} \sum_{k=0}^{\infty} \frac{y_{0,k}}{x^{2k}},
\end{gather*}
where
\begin{gather*}
y_{0,0}=A,\qquad y_{0,1}= \frac{\alpha A+\beta}{2}.
\end{gather*}

\item For each $j \in \{0,1,2,3\}$, there is a one-parameter family $\mathcal{F}_{A,j}$ of tronqu{\'e}e solutions of~\eqref{P4} in~$S^{(j)}$, where
\begin{gather*}
y(x) = l(x)+ x^{-1} h(w), \qquad w=x^2.
\end{gather*}
If $j$ is even, then $h(w)$ has the representations
\begin{gather*}
h(w) =
\begin{cases}
\displaystyle h^{+}(w) + \sum_{k=1}^{\infty}C_{+}^k {\rm e}^{-kw} w^{kM_1} \mathcal{L}_{\phi}{H}_k(w), & \displaystyle -\phi \in \left(0,\frac{\pi}{2}\right),\\
\displaystyle h^{-}(w) + \sum_{k=1}^{\infty}C_{-}^k {\rm e}^{-kw} w^{kM_1} \mathcal{L}_{\phi}{H}_k(w), & \displaystyle -\phi \in \left(-\frac{\pi}{2},0\right).
\end{cases}
\end{gather*}
If $j$ is odd, then $h(w)$ has the representations
\begin{gather*}
h(w) =
\begin{cases}
\displaystyle h^{-}(w) + \sum_{k=1}^{\infty}\hat{C}_{+}^k {\rm e}^{kw} (-w)^{kM_2} \mathcal{L}_{\phi}{\hat{H}}_k(-w), & \displaystyle -\phi \in \left(0,\frac{\pi}{2}\right),\\
\displaystyle h^{+}(w) + \sum_{k=1}^{\infty}\hat{C}_{-}^k {\rm e}^{kw} (-w)^{kM_2} \mathcal{L}_{\phi}{\hat{H}}_k(-w), & \displaystyle -\phi \in \left(-\frac{\pi}{2},0\right).
\end{cases}
\end{gather*}

\item Let $y(x)$ be a tronqu{\'e}e solution in $\mathcal{F}_{A,j}$. If $j$ is even, then the region of analyticity contains the corresponding branch of $(S_{\rm an}(h))^{1/2}$. If $j$ is odd, then the region of analyticity contains the corresponding branch of $\big(\hat{S}_{\rm an}(h)\big)^{1/2}$, and
\begin{gather*}
y(x) \sim \tilde{y}_0(x),\qquad x\in d \subset S^{(j)}, \qquad |x|\to\infty.
\end{gather*}
\end{enumerate}
\end{Theorem}

\begin{Theorem}\quad 
\begin{enumerate}[$(i)$]\itemsep=0pt
\item Let $j\in\{0,2\}$. For each $\delta>0$ there exists $R$ large enough such that we have a tritronqu{\'e}e solution $y^+_{j}$ analytic in $T^{(j)}_{\delta,R}$ given by
\begin{gather*}
y^+_{j}(x) = \frac{A}{ x}+ \frac{\alpha A+\beta }{2x^3} +h^{+}\big(x^2\big).
\end{gather*}

\item Let $j\in\{1,3\}$. For each $\delta>0$ there exists $R$ large enough such that we have a tritronqu{\'e}e solution $y^-_{j}$ analytic in $T^{(j)}_{\delta,R}$ given by
\begin{gather*}
y^-_{j}(x) = \frac{A}{ x}+ \frac{\alpha A+\beta }{2x^3} +h^{-}\big(x^2\big).
\end{gather*}
\end{enumerate}
\end{Theorem}

\begin{Note}In this case, the corresponding $F_0$ and $\hat{F}_0$ turn out to be $\xi$, which yield no singularities for $h$. However, it does not imply that the poles are nonexistent. More research needs to be done for this case.
\end{Note}

\section{Proofs and further results}
\subsection{Proof of Proposition~\ref{NFformal}}

Let $h$ and $\mathbf{u}$ be as defined in Section~\ref{intro}. We have a system of differential equations \eqref{NF} for $\mathbf{u}$. It is known (see \cite{OCIMRN,OCDuke, OCInv}) that it admits transseries solutions (i.e., formal exponential power series solutions) of the form
\begin{gather}\label{NFsys}
\tilde{\mathbf{u}}(w) = \tilde{\mathbf{u}}_0(w) +\sum_{k=1}^{\infty} C^k {\rm e}^{-kw}w^{-\beta_1k}\tilde{\mathbf{u}}_{k}(w),
\end{gather}
where $\tilde{\mathbf{u}}_0(w)$ and $\tilde{\mathbf{u}}_{k}(w)$ are formal power series in $w^{-1}$, namely
\begin{gather*}
\tilde{\mathbf{u}}_{k}(w) = \sum_{r=0}^{\infty}\frac{{\mathbf{u}}_{k,r}}{w^r},\qquad k \geq 1,\qquad
\tilde{\mathbf{u}}_{0}(w) = \sum_{r=2}^{\infty}\frac{{\mathbf{u}}_{0,r}}{w^r}.
\end{gather*}
Also, $\tilde{\mathbf{u}}_{0}(w)$ is the unique power series solution of~\eqref{NF}. The coefficients in the series $\tilde{\mathbf{u}}_{k}$ can be determined by substitution of the formal exponential power series $\tilde{\mathbf{u}}(w)$ into~\eqref{NFsys} and identification of each coefficient of~${\rm e}^{-kw}$. Proposition~\ref{NFformal} is then obtained through~\eqref{htou}. Furthermore,
\begin{gather}\label{serieshtou}
\tilde{h}_0(w) = \mathbf{r}_1 \cdot \tilde{\mathbf{u}}_{0}(w),\qquad
\tilde{s}_k(w) =\mathbf{r}_1 \cdot \tilde{\mathbf{u}}_{k}(w),\qquad \mathbf{r}_1 = \begin{bmatrix}
1-\dfrac{\beta_1}{2w}&1+\dfrac{\beta_2}{2w}\end{bmatrix}.
\end{gather}

\subsection{Proof of Theorem~\ref{NFactual}}

Let $d = {\rm e}^{{\rm i}\theta}\mathbb{R}^+$ with $\cos \theta>0$, and let $\mathbf{u}$ be a solution to~\eqref{NF} on $d$ for $w$ large enough, satisfying
\begin{gather*}
\mathbf{u}(w) \to \mathbf{0},\qquad w\in d,\qquad |w| \to \infty.
\end{gather*}
Theorem 3 in \cite{OCDuke}, Theorem 16, Lemma 17 and Theorem 19 in~\cite{OCInv} imply the following results:
\begin{Proposition}\label{sysu}\quad
\begin{enumerate}[$(i)$]\itemsep=0pt
\item For any $d'= {\rm e}^{{\rm i}\theta'}\mathbb{R}^+$ where $\cos \theta'>0$, the solution $\mathbf{u}(w)$ is analytic on $d'$ for $w$ large enough and $\mathbf{u} \sim \tilde{\mathbf{u}}_0(w)$ on~$d'$.
\item Given $\phi \in \big({-}\frac{\pi}{2},0\big) \cup\big(0, \frac{\pi}{2}\big)$, there exists a unique constant~$C(\phi)$ such that $\mathbf{u}$ has the following representation:
\begin{gather}\label{solnsys}
\mathbf{u}(w) = \mathcal{L}_{\phi} \mathbf{U}_0(w) + \sum_{k=1}^{\infty} (C(\phi) )^k {\rm e}^{-kw} w^{kM_1} \mathcal{L}_{\phi}\mathbf{U}_k(w),
\end{gather}
where
\begin{gather}\label{solnsys1}
\mathbf{U}_0 = \mathcal{B} \tilde{\mathbf{u}}_{0},\qquad
\mathbf{U}_k = \mathcal{B}\big(w^{-k\beta_1-kM_1}\tilde{\mathbf{u}}_{k}\big),\qquad k=1,2,\dots,
\end{gather}
where for each $k\geq 1$, $\mathbf{U}_k$ is analytic in the Riemann surface of $\mathbb{C}\backslash\left(\mathbb{Z}^+\cup \mathbb{Z}^-\right)$, and the branch cut for $\mathbf{U}_k$ is chosen to be $(-\infty,0]$. The function $C(\phi)$ is constant on $\big({-}\frac{\pi}{2},0\big)$ and also constant on $\big(0, \frac{\pi}{2}\big)$.
\item Let $\epsilon$ be small. There exist $\delta$, $R>0$ such that $\mathbf{u}(w)$ is analytic on
\begin{gather*}
S_{an,\epsilon}(\mathbf{u}(w)) = S^+_{\epsilon}\cup S^-_{\epsilon},
\end{gather*}
where $S^{\pm}$ is as defined in~\eqref{Spm}.
\end{enumerate}
\end{Proposition}

We now return to the proof of Theorem~\ref{NFactual}. Assume that $h(w)$ is a solution of~\eqref{eqh} on $d = {\rm e}^{{\rm i}\phi} \mathbb{R}^+$ with $\cos \phi>0$ for $|w|>w_0$, where $w_0>0$ is large enough. Without loss of generality we may assume that $w_0>\frac{\sqrt{|\beta_1 \beta_2|}}{2}$. Thus the vector function $\mathbf{u}(w)$ defined by
\begin{gather}\label{utoh}
\mathbf{u}(w) =\begin{bmatrix}
1-\dfrac{\beta_1}{2w}&1+\dfrac{\beta_2}{2w}\vspace{1mm}\\
-1-\dfrac{\beta_1}{2w}&1-\dfrac{\beta_2}{2w}
\end{bmatrix}^{-1}\begin{bmatrix}
h(w)\\
h'(w)
\end{bmatrix}
\end{gather}
is a solution of the differential system \eqref{NF}, and $h(w) = \mathbf{r}_1 \cdot \mathbf{u}(w)$.

Next we use the basic properties (see Lemmas \ref{L51} and~\ref{L52}) of the operators $\mathcal{B}$ and $\mathcal{L}_{\phi}$ and obtain the following
\begin{gather}
\mathcal{L}_{\phi}\mathcal{B} \big(\mathbf{r}_1\cdot \mathbf{\tilde{u}}_0\big) = \mathbf{r}_1\cdot \mathcal{L}_{\phi}\mathcal{B} \big( \mathbf{\tilde{u}}_0\big) =\mathbf{r}_1\cdot \mathcal{L}_{\phi} \mathbf{U}_0,\nonumber \\
\mathcal{L}_{\phi}\mathcal{B} \big(w^{-k\beta_1-kM_1}\mathbf{r}_1\cdot \mathbf{\tilde{u}}_k\big) = \mathbf{r}_1\cdot \mathcal{L}_{\phi}\mathcal{B} \big( w^{-k\beta_1-kM_1}\mathbf{\tilde{u}}_k\big)=\mathbf{r}_1\cdot \mathcal{L}_{\phi} \mathbf{U}_k.\label{r1pull}
\end{gather}
By Proposition~\ref{sysu}(i), given a ray $d'$ in the right half $w$-plane, $h(w) = \mathbf{r}_1 \cdot \mathbf{u}(w)$ is analytic on~$d'$ for $|w|$ large enough and is asymptotic to $\tilde{h}_0(w) = \mathbf{r}_1 \cdot \mathbf{\tilde{u}}_0(w)$ on $d'$. From the representations~\eqref{solnsys},~\eqref{solnsys1} in Proposition~\ref{sysu}(ii) of~$\mathbf{u}(w)$ and~\eqref{r1pull} we obtain the representations for $h(w)=\mathbf{r}_1 \cdot \mathbf{u}(w)$ as in~\eqref{truesolnup} and~\eqref{truesolndown}. For $|w|$ large enough, $h(w)$ is analytic where~$\mathbf{u}(w)$ is analytic, hence Proposition~\ref{sysu}(iii) implies Theorem~\ref{NFactual}(ii). Thus Theorem~\ref{NFactual} is proved.

\subsection{Proof of Theorem~\ref{sing}}

Let $h(w)$ and $\mathbf{u}(w)$ be as in the proof of Theorem~\ref{NFactual}. $h(w)$ has representations~\eqref{truesolnup} and~\eqref{truesolndown}. We will consider the case $C_+ \neq 0$ and prove~(i) and~(ii). The statements~(iii) and~(iv) about the case $C_- \neq 0$ follow by symmetry.

By Theorem~1 in~\cite{OCInv}, there exists $\delta_1>0$ such that for $|\xi_+|<\delta_1$ the power series
\begin{gather*}
\mathbf{G}_m(\xi_+) = \sum_{k=0}^{\infty}\xi_+^k\mathbf{u}_{k,m},\qquad m=0,1,2,\dots
\end{gather*}
converges, where $\mathbf{G}_0$ satisfies
\begin{gather*}
\mathbf{G}_0 =\mathbf{0},\qquad \mathbf{G}'_0 =\mathbf{e}_1.
\end{gather*}
Furthermore,
\begin{gather}\label{usim}
\mathbf{u}(w) \sim \sum_{m=0}^{\infty}w^{-m}\mathbf{G}_m\left(\xi_+(w)\right),\qquad |w| \to \infty
\end{gather}
holds uniformly in
\begin{gather*}
\mathcal{S}_{\delta_1} = \left\{w\colon \arg(w) \in \left(-\frac{\pi}{2}+\delta,\frac{\pi}{2}+\delta \right), \, |\xi_+(w)|<\delta_1\right\}.
\end{gather*}
By Theorem 2 in \cite{OCInv}, for $R$ large enough and $\delta$, $\epsilon$ small enough, $\mathbf{u}(w)$ is analytic in $\mathcal{D}_w^+$ (see~\eqref{Dwp}). Also, the asymptotic representation~\eqref{usim} holds in~$\mathcal{D}_w^+$. Moreover, if~$\mathbf{G}_0$ has an isolated singularity at $\xi_s$, then $\mathbf{u}(w)$ is singular at a distance at most $o(1)$ of $w_n^+$ given in~\eqref{wpn}, as $w_n^+ \to \infty$.
Since $h(w) = \mathbf{r}_1 \cdot \mathbf{u}(w)$, Theorem~\ref{sing}(i) follows from the results cited.

Assume $|w|> {\sqrt{|\beta_1 \beta_2|}}\slash{2}$. Both \eqref{htou} and \eqref{utoh} hold. While \eqref{htou} implies that $h$ is analytic at least where $\mathbf{u}$ is analytic, \eqref{utoh} implies that $h$ is singular where $\mathbf{u}$ is singular. Thus the asymptotic position of singularities, i.e., poles of $h(w)$ is the same as that of $\mathbf{u}(w)$, which is presented in equation~\eqref{wpn}. Thus Theorem~\ref{sing}(ii) is proved.

\subsection{Proof of Corollary~\ref{NFtri}}

Let the notations be the same as in Section~\ref{NFtrit}. First we point out some properties of $\mathbf{U}_0(p)$. See also~\cite{OCIMRN}.

We apply formal inverse Laplace transform to the system \eqref{NF}. To be precise, assume the analytic function $\mathbf{g}(w,\mathbf{u})$ has the Taylor expansion at $(\infty,\mathbf{0})$ as follows
\begin{gather*}
\mathbf{g}(w,\mathbf{u}) = \sum_{m\geq0;|\mathbf{l}|\geq 0} \mathbf{g}_{m,\mathbf{l}}w^{-m} \mathbf{u}^{\mathbf{l}}, \qquad \big\vert w^{-1}\big\vert < \xi_0, \qquad \vert \mathbf{u} \vert < \xi_0.
\end{gather*}
Note that by assumption $\mathbf{g}_{m,\mathbf{l}}=\mathbf{0}$ if~$|\mathbf{l}|\leq 1$ and $m\leq 1$. Denote $\mathbf{U} = \mathcal{L}^{-1}\mathbf{u}$. Then the formal inverse Laplace transform of the differential system~\eqref{NF} is the system of convolution equations
\begin{gather}\label{ILTNF}
-p \mathbf{U}(p) = -\left[\hat{\Lambda} \mathbf{U}(p)+\hat{B} \int_{0}^{p}\mathbf{U}(s)\mathrm{d}s\right] +\mathcal{N}\left(\mathbf{U}\right)(p),
\end{gather}
where
\begin{gather*}
\mathcal{N}\left(\mathbf{U}\right)(p) = \sum_{m=2}^{\infty} \frac{\mathbf{g}_{m,\mathbf{0}}}{(m-1)!}p^{m-1}+\sum_{|\mathbf{l}|\geq 2} \mathbf{g}_{0,\mathbf{l}} \mathbf{U}^{\ast \mathbf{l}}+\sum_{|\mathbf{l}|\geq 1}\left(\sum_{m=1}^{\infty}\frac{\mathbf{g}_{m,\mathbf{l}}}{(m-1)!}p^{m-1}\right) \ast \mathbf{U}^{\ast \mathbf{l}}.
 \end{gather*}
Let $\mathbf{v}(p) = \left(v_1(p),\dots,v_n(p)\right)$ be an $n$-dimensional complex vector function, $f(p)$~be a locally integrable complex function and $\mathbf{l} = (l_1,\dots,l_n)$ be an $n$-dimensional multi-index. Then
\begin{gather*}
\mathbf{v}^{\ast \mathbf{l}} := v_1^{\ast l_1} \ast v_2^{\ast l_2} \ast\cdots \ast v_n^{\ast l_n},\qquad
\left(\mathbf{v}\ast f\right) (p) \in \mathbb{C}^{n}, \qquad (\mathbf{v} \ast f )_{j} = v_j \ast f, \qquad j=1,\dots,n.
\end{gather*}
We gather the following facts about $\mathbf{U}_0$.

\begin{Proposition}\label{U0}\quad
\begin{enumerate}[$(i)$]\itemsep=0pt
\item Let ${K} \in \mathcal{O}$ be a closed set such that for every point $p \in K$, the line segment connecting the origin and $p$ is contained in~$K$. Then $\mathbf{U}_0$ is the unique solution to \eqref{ILTNF} in $K$.
\item $\mathbf{U}_0 = \mathcal{B} \tilde{\mathbf{u}}_0$. $\mathbf{U}_0$ is analytic in the domain $\mathcal{O}= \mathbb{C} \backslash \left[(\infty,-1] \cup [1,\infty)\right]$, and is Laplace transformable along any ray ${\rm e}^{{\rm i}\phi} \mathbb{R}^+$ contained in $\mathcal{O}$. $\mathcal{L}_{\phi} \mathbf{U}_0$ is a solution of~\eqref{ILTNF} for each~$\phi$ such that $|\cos(\phi)|<1$.
\item
Let $K$ be as in $(i)$. There exists $b_K>0$ large enough such that
\begin{gather*}
\sup_{p \in K} \int_{[0,p]} \vert \mathbf{U}_0(s) \vert {\rm e}^{-b_K |s|} |\mathrm{d}s| < \infty.
\end{gather*}
\end{enumerate}
\end{Proposition}

Proposition~\ref{U0}(i) and (ii) come from Proposition 6 in~\cite{OCIMRN}. Although (iii) is not stated expli\-cit\-ly in~\cite{OCIMRN}, it can be easily obtained by the same approach used to prove Proposition~6. Let $K$ be as in~(i). Consider the Banach space
\begin{gather*}
L_{\rm ray}(K):=\left\{\mathbf{f}\colon \mathbf{f} \ \text{is locally integrable on }[0,p]\text{ for each }p\in K \right\}
\end{gather*}
equipped with the norm $\| \cdot \|_{b,K}$ defined by
\begin{gather*}
\|\mathbf{f}\|_{b,K}:=\sup_{p\in K} \int_{[0,p]}\| \mathbf{f}(s) \| {\rm e}^{-b |s|} |\mathrm{d}s|,
\end{gather*}
where $\| \mathbf{f}(s) \| = \max\{|f_1(s)|,|f_2(s)|\}$. We can show that for $b$ large enough, the operator
\begin{gather*}
\mathcal{N}_1:=\mathbf{U}(p) \mapsto \big(\hat{\Lambda} - p I\big)^{-1} \left(-\hat{B} \int_{0}^{p}\mathbf{U}(s) \mathrm{d}s+\mathcal{N} (\mathbf{U} )(p)\right)
\end{gather*}
is contractive in the closed ball $\mathcal{S}:= \{\mathbf{f}\in L_{\rm ray}(K)\colon \|\mathbf{f}\|_{b,K} \leq \delta\}$ of $L_{\rm ray}(K)$ if $\delta$ is small enough. By contractive mapping theorem there is a unique solution of $\mathcal{N}_1 \mathbf{U} =\mathbf{U} $ in $\mathcal{S}$, namely~$\mathbf{U}_0$ by uniqueness of the solution. Using integration by parts and (iii) we have the following:

\begin{Corollary}\label{LU0facts}\quad
\begin{enumerate}[$(i)$]\itemsep=0pt
\item If $\phi \in (0,\pi)$ or $\phi \in (-\pi,0)$, $\mathcal{L}_{\phi} \mathbf{U}_0(w)$ is analytic $($at least$)$ in the region
\begin{gather*}
\mathcal{A}_{\phi}:= \{w\colon |w|\cos (\phi+\arg(w) )>b \},
\end{gather*}
where $b=b_K$ is as in Proposition~{\rm \ref{U0}}$(iii)$ with $K={\rm e}^{{\rm i}\phi}\mathbb{R}^+$.
\item If $0<\phi_1< \phi_2 < \pi$ or $0<-\phi_1< -\phi_2 < \pi$, then $\mathcal{L}_{\phi_1}\mathbf{U}_0$ and $\mathcal{L}_{\phi_2}\mathbf{U}_0$ are analytic continuations of each other.
\end{enumerate}
\end{Corollary}

Since $H_0 = \mathcal{B} \tilde{h}_0$, by \eqref{serieshtou} and Lemma~\ref{L51} we have
\begin{gather}
H_0(p) = (\mathcal{B} u_{0,1} )(p)+ (\mathcal{B} u_{0,2} )(p) -\frac{\beta_1}{2} [1\ast (\mathcal{B} u_{0,1} ) ](p)+\frac{\beta_2}{2} [1\ast (\mathcal{B} u_{0,2} ) ](p)\nonumber\\
\hphantom{H_0(p)}{} =U_{0,1}(p)+U_{0,2}(p) -\frac{\beta_1}{2} (1\ast U_{0,1} )(p)+\frac{\beta_2}{2} (1\ast U_{0,2} )(p),\label{H0U0}
\end{gather}
where ${u}_{0,i}$, $i = 1,2$, is the $i$-th component of the vector function $\mathbf{u}_0$ and ${U}_{0,i}$, $i = 1,2$, is the $i$-th component of $\mathbf{U}_0$. It is clear from \eqref{H0U0} that Proposition~\ref{U0}(iii) and Corollary~\ref{LU0facts} hold with~$\mathbf{U}_0$ replaced by~$H_0$.

Merely by $\hat{H}_0(p) = -H_0(-p)$ and Corollary~\ref{LU0facts}(ii) with $\mathbf{U}_0$ replaced by $H_0$ we obtain Corollary~\ref{NFtri}(i). Moreover, both~$h^+$ and~$h^-$ are special cases of tronqu{\'e}e solutions, thus Theorems~\ref{NFactual} and~\ref{sing} apply. $h^+$ is analytic at least on $S_{\rm an}(h^+)\cup\big({-}S_{\rm an}\big(\hat{h}^+\big)\big) $ and $h^-$ is analytic at least on $S_{\rm an}(h^-)\cup\big({-}S_{\rm an}\big(\hat{h}^-\big)\big) $. We also obtain the asymptotic position of singularities of the tritronqu{\'e}e solutions as in Corollary~\ref{NFtri}(ii).

\subsection{Proof of the results in Section~\ref{BS3}}

Once we have the normalizations in the form of \eqref{eqh} of the equations \eqref{P31}, \eqref{P32} and \eqref{P4}, the results in Section~\ref{BS3} follow from the results in Section~\ref{NFinfo}. Here we present the details of finding solutions to \eqref{P31eqF0}, \eqref{P32eqF0}, \eqref{P41eqF0} and \eqref{P42eqF0} satisfying \eqref{conF0}.

\subsubsection{Solving (\ref{P31eqF0})}
Make the substitution $Q(s) = A+F_0({\rm e}^s)$ then \eqref{P31eqF0} transforms into
\begin{gather*}
{\frac {\mathrm{d}^{2}}{\mathrm{d}{s}^{2}}}Q (s) -{\frac { \big( {\frac {\mathrm{d}
}{\mathrm{d}s}}Q (s) \big) ^{2}}{Q (s) }}- {\frac { (Q(s)) ^{3}}{4{A}^{2}}}+ {\frac {1}{4{A}^{2}Q (s) }} =0.
\end{gather*}
Multiplying both sides by $1/Q(s)$ we obtain
\begin{gather*}
\frac{\mathrm{d}}{\mathrm{d}s}\left(\frac{Q'(s)}{Q(s)}\right) = \frac{1}{4A^2}\left(Q^2(s)-\frac{1}{Q^2(s)}\right).
\end{gather*}
Multiplying both sides by ${2Q'(s)}/{Q(s)}$ and integrating with respect to $s$ we have
\begin{gather}
\left(\frac{Q'(s)}{Q(s)}\right)^2 = \frac{1}{4A^2}\left(Q^2(s)+\frac{1}{Q^2(s)}+C_1\right),\qquad \text{i.e.},\nonumber\\
(Q'(s))^2 = \frac{1}{4A^2}\big(Q^4(s)+C_1Q^2(s)+1\big).\label{P31eqQ}
\end{gather}
By a linear transformation $Q(s) = \tilde{Q}(s)/(2A)$, \eqref{P31eqQ} is reduced to the Jacobi normal form which is solved by Jacobi elliptic functions unless $C_1 \in \{-2,2\}$. Since $Q(s) = A+F_0({\rm e}^s)$ and $F_0(0)=0$ (see~\eqref{conF0}), the solution we look for cannot be an elliptic function. Moreover, as $\operatorname{Re}(s) \to -\infty$, $Q(s) \to A$ implies $Q'(s) \to 0$, so~\eqref{P31eqQ} needs to be of the form
\begin{gather*}
(Q'(s))^2 = \frac{1}{4A^2}\big(Q^2(s)-A^2\big)^2.
\end{gather*}
The solution satisfying $Q(s) \to A$ as $\operatorname{Re}(s) \to -\infty$ is
\begin{gather*}
Q(s) = A \cdot\frac{C_2-{\rm e}^s}{C_2+{\rm e}^s}, \qquad C_2\neq 0.
\end{gather*}
Thus the solution to \eqref{P31eqF0} is
\begin{gather*}
F_0(\xi) =-\frac{2A\xi}{C_2+\xi}.
\end{gather*}
Hence the solution to~\eqref{P31eqF0} satisfying~\eqref{conF0} is
\begin{gather*}
F_0(\xi) =\frac{2A\xi}{2A-\xi}.
\end{gather*}

\subsubsection{Solving (\ref{P32eqF0})}
Make the substitution $Q(s) = A+F_0({\rm e}^s)$ then \eqref{P32eqF0} transforms into
\begin{gather*}
{\frac {{\rm d}^{2}}{{\rm d}{s}^{2}}}Q (s) -{\frac {\big( {\frac {\rm d}{{\rm d}s}}Q (s) \big) ^{2}}{Q
(s) }}- {\frac { (Q(s))^{2}}{3A}}+ {\frac {1}{3AQ (s) }}=0.
\end{gather*}
Multiplying both sides by $1/Q(s)$ we obtain
\begin{gather*}
\frac{\mathrm{d}}{\mathrm{d}s}\left(\frac{Q'(s)}{Q(s)}\right) = \frac{1}{3A}\left(Q(s)-\frac{1}{Q^2(s)}\right).
\end{gather*}
Multiplying both sides by ${2Q'(s)}/{Q(s)}$ and integrating with respect to $s$ we have
\begin{gather}
\left(\frac{Q'(s)}{Q(s)}\right)^2 = \frac{1}{3A}\left(2Q(s)+\frac{1}{Q^2(s)}+C_1\right),\qquad\text{i.e.},\nonumber\\
(Q'(s))^2 = \frac{1}{3A}\big(2Q^3(s)+C_1Q^2(s)+1\big).\label{P32eqQ}
\end{gather}
Notice that if the equation $2x^3+C_1x^2+1 =0$ has three distinct roots then~\eqref{P32eqQ} is known to have Weierstrass $\wp$-functions as general solutions, in which case the corresponding $F_0(\xi) = Q(\ln \xi)-A$ fails to satisfy the condition~\eqref{conF0}. Hence $C_1$ must be such that the equation $2x^3+C_1x^2+1 =0$ has a multiple root. Denote the multiple root by~$r_1$. Then
\begin{gather*}
2x^3+C_1x^2+1 =2(x-r_1)^2(x-r_2).
\end{gather*}
Then we obtain
\begin{gather*}
r_1 = -2 r_2, \qquad r^3_1=1, \qquad C_1 = -3 r_1.
\end{gather*}
Since $Q(s) \to A$ as $\operatorname{Re}(s) \to -\infty$, $Q'(s) \to 0$. Hence $r_1 = A$ or $r_2 = A$. We knew from the normalization (see~\eqref{P32lx}) that $A^3=1$. Thus $r_1=A$, $r_2 = -A/2$, and $C_1 = -3A$. Hence~\eqref{P32eqQ} is of the form
\begin{gather}\label{P32eqQ2}
(Q'(s))^2 = \frac{2}{3A} (Q(s)-A )^2\left(Q(s)+\frac{A}{2}\right).
\end{gather}
The solution to \eqref{P32eqQ2} is
\begin{gather*}
Q(s) = -\frac{A}{2}+\frac{3A}{2}\left(\frac{C_2-{\rm e}^s}{C_2+{\rm e}^s}\right)^2.
\end{gather*}
Hence the solution of \eqref{P32eqF0} satisfying \eqref{conF0} is
\begin{gather*}
F_0(\xi) = \frac{36A^2\xi}{(6A-\xi)^2}.
\end{gather*}

\subsubsection{Solving (\ref{P41eqF0})}
Make the substitution $Q(s) = F_0({\rm e}^s)-2/3$ then \eqref{P41eqF0} transforms into
\begin{gather*}
{\frac {{\rm d}^{2}}{{\rm d}{s}^{2}}}Q (s) -{\frac {\big( {\frac {\rm d}{{\rm d}s}}Q (s) \big)^2}{2 Q
(s) }}+ {\frac { 9 Q^3 (s) }{8}}+3 Q^2(s)+ {\frac {3 Q(s)}{2}}=0.
\end{gather*}
Multiplying both sides by $2Q'(s)/Q(s)$ and we have
\begin{gather*}
\frac{{\rm d}}{{\rm d}{s}}\left[\frac{(Q'(s))^2}{Q(s)}\right] = \frac{{\rm d}}{{\rm d}{s}}\left(-\frac{3}{4}Q^3(s)-3 Q^2(s)-3 Q(s)\right).
\end{gather*}
Integrating with respect to $s$ we have
\begin{gather}\label{P41eqQ}
(Q'(s))^2 = -\frac{3}{4}Q^4(s)-3 Q^3(s)-3 Q^2(s) +C_1 Q(s).
\end{gather}
Letting $\operatorname{Re}(s) \to -\infty$ we have $Q(s) \to -2/3$ and $Q'(s) \to 0$. Thus $C_1 = -8/9$ and the equation~\eqref{P41eqQ} is of the form
\begin{gather*}
(Q'(s))^2 = -\frac{3}{4} Q(s)\left(Q(s)+\frac{8}{3}\right)\left(Q(s)+\frac{2}{3}\right)^2.
\end{gather*}
This is a separable differential equation with general solutions
\begin{gather*}
Q(s) = -\frac{2}{3} \frac{{\rm e}^{2s}-C_2^2-2{\rm i}C_2 {\rm e}^{s}}{{\rm e}^{2s}-C_2^2+{\rm i}C_2{\rm e}^{s}}.
\end{gather*}
Hence the solution of \eqref{P41eqF0} satisfying \eqref{conF0} is
\begin{gather*}
F_0(\xi) = \frac{4 \xi}{\xi^2+2\xi+4}.
\end{gather*}

\subsubsection{Solving (\ref{P42eqF0})}
Make the substitution $Q(s) = F_0({\rm e}^s)-2$ then \eqref{P41eqF0} transforms into
\begin{gather*}
{\frac {{\rm d}^{2}}{{\rm d}{s}^{2}}}Q (s) -{\frac {\big( {\frac {\rm d}{{\rm d}s}}Q (s) \big)^2}{2 Q(s) }}- {\frac { 3 Q^3 (s) }{8}}- Q^2(s)- {\frac {Q(s)}{2}}=0.
\end{gather*}
Multiplying both sides by $2Q'(s)/Q(s)$ and we have
\begin{gather*}
\frac{{\rm d}}{{\rm d}{s}}\left[\frac{(Q'(s))^2}{Q(s)}\right] = \frac{{\rm d}}{{\rm d}{s}}\left(\frac{1}{4}Q^3(s)+ Q^2(s)+ Q(s)\right).
\end{gather*}
Integrating with respect to $s$ we have
\begin{gather}\label{P41eqQn}
(Q'(s))^2 = \frac{1}{4}Q^4(s)+ Q^3(s)+ Q^2(s) +C_1 Q(s).
\end{gather}
Letting $\operatorname{Re}(s) \to -\infty$ we have $Q(s) \to -2$ and $Q'(s) \to 0$. Thus $C_1 = 0$ and the equation~\eqref{P41eqQn} is of the form
\begin{gather*}
(Q'(s))^2 = \frac{1}{4} Q^2(s)(Q(s)+2)^2.
\end{gather*}
This differential equation has general solutions
\begin{gather*}
Q(s) = -\frac{2 C_2}{C_2+{\rm e}^{s}}.
\end{gather*}
Hence the solution of \eqref{P42eqF0} satisfying \eqref{conF0} is
\begin{gather*}
F_0(\xi) = \frac{2 \xi}{\xi+2}.
\end{gather*}

\appendix

\section{Appendix}\label{Appx}

Recall that the Borel transform of a formal series
\begin{gather*}
\tilde{f}(w) = \sum_{n=0}^{\infty}a_n w^{-r-n},\qquad \operatorname{Re}(r)>0,
\end{gather*}
where the series $\sum\limits_{n=0}^{\infty}a_n x^{n}$ has a positive radius of convergence, is defined to be the formal power series
\begin{gather*}
\big(\mathcal{B} \tilde{f} \big)(p) := \sum_{n=0}^{\infty}\frac{a_n p^{n+r-1}}{\Gamma(n+r)}.
\end{gather*}

\begin{Lemma}\label{L51}Assume that we have two formal series $\tilde{f}$ and $\tilde{g}$,
\begin{gather*}
\tilde{f}(w) = \sum_{n=0}^{\infty}a_n w^{-r-n},\qquad \operatorname{Re}(r)>0,\\
\tilde{g}(w) = \sum_{n=0}^{\infty}b_n w^{-r-s},\qquad \operatorname{Re}(s)>0,
\end{gather*}
where both series $\sum\limits_{n=0}^{\infty}a_n x^{n}$ and $\sum\limits_{n=0}^{\infty}b_n x^{n}$ have positive radii of convergence. Then
\begin{gather*}
\mathcal{B}\big(\tilde{f} \tilde{g}\big)(p) = \big(\mathcal{B}\tilde{f} \ast \mathcal{B}\tilde{g}\big) (p)= p^{r+s-1} \sum_{n=0}^{\infty}\left(\sum_{k=0}^{n} a_k b_{n-k}\right)\frac{p^n}{\Gamma(n+r+s)},
\end{gather*}
where
\begin{gather*}
\big(\mathcal{B}\tilde{f} \ast \mathcal{B}\tilde{g}\big)(p) := \int_{0}^{p}\big(\mathcal{B}\tilde{f}\big)(t)\big(\mathcal{B}\tilde{g}\big)(p-t) \mathrm{d}t.
\end{gather*}
\end{Lemma}

Recall that the Laplace transform $\mathcal{L}_{\phi}$ is defined as the following
\begin{gather*}
f\longmapsto \int_{0}^{\infty {\rm e}^{{\rm i}\phi}} f(p){\rm e}^{-xp}{\rm d}p,
\end{gather*}
where $\phi \in \mathbb{R}$.

\begin{Lemma}\label{L52}
Assume that the function $f$ is integrable over the ray ${\rm e}^{{\rm i}\phi} \mathbb{R}^+$, namely
\begin{gather*}
\int_{0}^{\infty {\rm e}^{{\rm i}\phi} } \vert f(p) \vert \vert \mathrm{d}p \vert < \infty.
\end{gather*}
Then for $\operatorname{Re}(w {\rm e}^{{\rm i}\phi})>0 $,
\begin{gather*}
\mathcal{L}_{\phi}(1\ast f) (w) = \frac{1}{w} \mathcal{L}_{\phi}(f) (w).
\end{gather*}
\end{Lemma}

\subsection*{Acknowledgements}

I am very grateful for the advice and help of Professors Ovidiu Costin and Rodica Costin when I worked on this problem. I also greatly appreciate the referees whose comments helped me improve my paper significantly.

\pdfbookmark[1]{References}{ref}
\LastPageEnding

\end{document}